\documentclass[a4paper,twoside,10pt]{article}
\pdfoutput=1

	\usepackage{amsmath}				
	\usepackage{amscd}					
	\usepackage{amsthm,thmtools}		
	\usepackage{amssymb}				
	\usepackage[titletoc]{appendix}		
	\usepackage{titling}				
	\usepackage[
		style=alphabetic,				
		sorting=nyt,					
		maxnames=5,						
		minnames=3,						
		backend=bibtex8,				
		firstinits=true,				
		doi=false,						
	]{biblatex}
	\usepackage{bm}						
	\usepackage{color}					
	\usepackage{enumitem}				
	\usepackage{etoolbox}				
	\usepackage{fancyhdr}				
	\usepackage{float}					
	\usepackage[T1]{fontenc}			
	\usepackage{gensymb}				
	\usepackage[a4paper,
		outer=3cm,
		inner=3cm,
		top=3cm,
		bottom=3cm]{geometry}
	\usepackage[utf8]{inputenc}			
	\usepackage{ifmtarg}				
	\usepackage{mathrsfs}				
	\usepackage{makeidx}				
	\usepackage{multicol}				
	\usepackage{parskip}				
	\usepackage[refpage,intoc]{nomencl}		
	\usepackage{slashed}				
	\usepackage{stmaryrd}				
	\usepackage{tensor}					
	\usepackage{textcomp}				
	\usepackage{titlesec}               
	\usepackage{tikz}                   
	\usepackage{tocloft}				
	\usepackage{url}					
	\usepackage[all,cmtip]{xy}			
	
	\usepackage{csquotes}				
	
	\usepackage[pageanchor,breaklinks]{hyperref}		
	\usepackage{cleveref}	

	\bibliography{refs} 
	\DeclareFieldFormat{postnote}{#1}
	\DeclareFieldFormat{multipostnote}{#1}
	
		\numberwithin{equation}{section}
		\crefformat{equation}{(#2#1#3)} 

	\fancypagestyle{plain}{%
		\fancyhf{}%
		\fancyfoot[C]{}%
	}

		\let\oldtitle\title
		\renewcommand{\title}[1]{\oldtitle{#1}\newcommand{\mythetitle}{#1}}
		\let\oldauthor\author
		\renewcommand{\author}[1]{\oldauthor{#1}\newcommand{\mytheauthor}{#1}}	
		\definecolor{linkcolor}{RGB}{00,10,138} 
		\hypersetup{
			colorlinks=true,
			breaklinks=true,
			linkcolor=linkcolor,
			menucolor=linkcolor,
			citecolor=linkcolor,
			filecolor=linkcolor,
			urlcolor=linkcolor,
			frenchlinks=false
		}
	
		\renewcommand{\nomname}{List of symbols}
		\patchcmd{\thenomenclature}
		  {\chapter*{\nomname}}
		  {\phantomsection\chapter*{\nomname}\label{SectNomenclature}}
		  {}
		  {}

		\usetikzlibrary{decorations.markings}

		\titleformat{\section}[block]{\large\bfseries\scshape\filcenter\textsc}{\thesection. }{0pt}{}
		\titleformat{\subsection}{\bfseries\scshape}{\thesubsection. }{0pt}{}
		
	\declaretheoremstyle[
		spaceabove=6pt, spacebelow=6pt,
		headfont=\normalfont\bfseries,
		notefont=\mdseries, notebraces={(}{)},
		bodyfont=\normalfont,
		postheadspace=1em,
		qed=$\lozenge$,
	]{mythmstyle}

	\declaretheoremstyle[
		spaceabove=6pt, spacebelow=6pt,
		headfont=\normalfont\bfseries,
		notefont=\mdseries, notebraces={(}{)},
		bodyfont=\normalfont,
		postheadspace=1em,
		qed=$\square$,
	]{myprfstyle}

	\declaretheorem[style=mythmstyle,name=Main Theorem]{MainThm}

	\declaretheorem[style=mythmstyle,name=Definition,numberwithin=section]{Def}
	
	\declaretheorem[style=mythmstyle,sibling=Def,name=Theorem]{Thm}

	\declaretheorem[style=mythmstyle,sibling=Def,name=Corollary]{Cor}
	
	\declaretheorem[style=mythmstyle,sibling=Def,name=Lemma]{Lem}
	
	\declaretheorem[style=mythmstyle,sibling=Def,name=Remark]{Rem}
	\declaretheorem[style=mythmstyle,sibling=Def,name=Conjecture]{Cnj}

	\declaretheorem[style=myprfstyle,numbered=no,name=Proof]{Prf}

	\crefname{Def}{Definition}{Definitions}
	\crefname{Exm}{Example}{Examples}
	\crefname{Con}{Convention}{Conventions}
	\crefname{Cor}{Corollary}{Corollaries}
	\crefname{Fct}{Fact}{Facts}
	\crefname{Lem}{Lemma}{Lemmas}
	\crefname{NumText}{Paragraph}{Paragraphs}
	\crefname{Rem}{Remark}{Remarks}
	\crefname{Cnj}{Conjecture}{Conjectures}
	\crefname{Exc}{Exercise}{Exercises}
	\crefname{MainThm}{Main Theorem}{Main Theorems}
	\crefname{Prb}{Problem}{Problems}
	\crefname{Qst}{Question}{Questions}
	\crefname{Thm}{Theorem}{Theorems}
	\crefname{Not}{Notation}{Notation}
	\crefname{Prf}{Proof}{Proofs}
	\crefname{enumi}{}{}

		\sbox0{$\mathbf{\xdef\mybffam{\the\fam}}$}

		\newlist{category}{enumerate}{10}
		\setlist[category]{
			itemsep=1ex,
			topsep=0.5ex,
			parsep=0pt,
			leftmargin=3em,
			itemindent=5em,
			labelwidth=5em,
			labelsep=1em,
			align=categorylabel,
			label*=.\arabic*
		}
		\setlist[category,1]{label=\arabic*}
		\SetLabelAlign{categorylabel}{\bfseries{\scshape{#1}}:\hfil}	
	
		\newenvironment{DefI}[1][]{\ifstrempty{#1}{\Def}{\Def[#1]\index{#1}}}{\endDef}
	
		\newcommand{\DefMap}[4]{
			\begin{align*}
				\begin{array}{rcl}
					#1 & \to & #2 \\
					#3 & \mapsto & #4 
				\end{array} 
			\end{align*}
		}

	\newlist{caselist}{enumerate}{10}
	\setlist[caselist]{itemsep=1ex,topsep=1ex,parsep=0pt,leftmargin=0pt,%
	  labelwidth=*,labelsep=1ex,align=caselabel,label*=.\arabic*}
	\setlist[caselist,1]{label=\arabic*}

	\newcommand\thecaselabeltext{}
	\SetLabelAlign{caselabel}{\textsc{Case~#1}\thecaselabeltext:\hfil}

	\crefname{caselisti}{case}{cases}
	\loop\ifnum\count255<10
	   \advance\count255 1
	   \crefalias{caselist\romannumeral\count255}{caselisti}
	\repeat

	\makeatletter
	\def\ifempty#1{\def\@x{#1}\ifx\@x\@empty}
	\makeatother

	\newlist{steplist}{enumerate}{10}
	\setlist[steplist]{nolistsep,itemsep=0ex,topsep=0ex,partopsep=0px,parsep=0pt,leftmargin=0pt,%
	  labelwidth=*,labelsep=1ex,align=steplabel,label*=.\arabic*}
	\setlist[steplist,1]{label=\arabic*}

	\newcommand\thesteplabeltext{}
	\SetLabelAlign{steplabel}{\textsc{Step~#1}\thesteplabeltext:\hfil}

	\crefname{steplisti}{Step}{steplist}
	\loop\ifnum\count255<10
	   \advance\count255 1
	   \crefalias{steplist\romannumeral\count255}{steplisti}
	\repeat

	\makeatletter
	\def\ifempty#1{\def\@x{#1}\ifx\@x\@empty}
	\makeatother

	\newcommand\step[1][]{\renewcommand{\thesteplabeltext}
	  {\ifempty{#1}\else~(#1)\fi}\item}
	\DeclareMathOperator{\C}{\mathbb{C}}			
	\DeclareMathOperator{\Diff}{Diff} 				
	\DeclareMathOperator{\Dirac}{\slashed{D}}		
	\newcommand{\GLp}{\operatorname{GL}^{+}}		
	\newcommand{\GLtp}{\widetilde{\operatorname{GL}}^{+}}	
	\DeclareMathOperator{\GL}{GL} 					
	\DeclareMathOperator{\I}{\operatorname{I}} 		
	\DeclareMathOperator{\K}{\mathbb{K}}			
	\DeclareMathOperator{\Lin}{span} 				
	\DeclareMathOperator{\N}{\mathbb{N}}			
	\DeclareMathOperator{\Rm}{\mathcal{R}}			
	\DeclareMathOperator{\R}{\mathbb{R}}			
	\DeclareMathOperator{\SO}{SO} 					
	\DeclareMathOperator{\Spin}{Spin} 				
	\DeclareMathOperator{\Z}{\mathbb{Z}}			
	\DeclareMathOperator{\dist}{dist} 				
	\DeclareMathOperator{\dom}{\mathcal{D}} 		
	\DeclareMathOperator{\id}{id} 					
	\DeclareMathOperator{\pr}{pr} 					
	\DeclareMathOperator{\sgn}{sgn} 				
	\DeclareMathOperator{\spec}{spec} 				
	\DeclareMathOperator{\spin}{spin} 				
	
	\newcommand{\subseto}{\hspace{.1em}\mathring{\subset}\hspace{.1em}} 
	\renewcommand{\H}{\mathbb{H}}			
	\renewcommand{\S}{S}					
	\newcommand{\D}{D}						

	\title{Existence of Dirac Eigenvalues of higher Multiplicity}
	\author{Nikolai Nowaczyk}
	
	\makeindex
	\usepackage[totoc]{idxlayout}
	
	\makenomenclature

	\nomenclature[N]{$\N$}{the natural numbers $\N=\{0,1,2, \ldots \}$}
	\nomenclature[TN]{$\mathcal{T}(M)$}{smooth vector fields on $M$}

\begin{document}

		\begin{center}
			{\Large \bfseries \mythetitle} \\
			\vspace{0.3cm}
			\mytheauthor \footnote{
				\textbf{Contact Details:} \\
				\textbf{Address:} Universität Regensburg, Fakultät für Mathematik, Universitätsstr. 31, 93040 Regensburg, Germany \\
				\textbf{Email:} mail@nikno.de} \\
			\vspace{0.3cm}
			\today \\
		\end{center}
		
	
		\vspace{1em}
		\textbf{Abstract.} In this article, we prove that on any compact spin manifold of dimension $m \equiv 0,6,7 \mod 8$, there exists a metric, for which the associated Dirac operator has at least one eigenvalue of multiplicity at least two. We prove this by ``catching'' the desired metric in a subspace of Riemannian metrics with a loop that is not homotopically trivial. We show how this can be done on the sphere with a loop of metrics induced by a family of rotations. Finally, we transport this loop to an arbitrary manifold (of suitable dimension) by extending some known results about surgery theory on spin manifolds. 
		
		\vspace{1em}
		
		\textbf{Keywords.} spin geometry, Dirac operator, spectral geometry,  Dirac spectrum, prescribing eigenvalues, surgery theory  
		
		\vspace{1em}
		
		\textbf{Mathematics Subject Classification 2010.} 53C27, 58J05, 58J50, 57R65
		
		\vspace{1em}

		\tableofcontents
		\thispagestyle{empty}
		

		\newpage
		
		\section{Introduction and statement of the results}
\label{SectHigherIntro}

\nomenclature[M]{$M$}{a closed spin manifold of dimension $m$}
\nomenclature[Theta]{$\Theta$}{a topological spin structure}
\nomenclature[SigmagM]{$\Sigma^g_{\K} M$}{spinor bundle over $M$ w.r.t. $g$}
\nomenclature[Kbb]{$\K$}{$\R$ or $\C$}
\nomenclature[DiracgK]{$\Dirac^g_{\K}$}{Dirac operator w.r.t. $g$ over $\K$}
\nomenclature[spec]{$\spec \Dirac^g_{\K}$}{Dirac spectrum}
\nomenclature[GammaSigma]{$\Gamma(\Sigma^g_{\K} M)$}{smooth spinor fields}
\nomenclature[L2Sigma]{$L^2(\Sigma^g_{\K} M)$}{$L^2$ spinor fields}
\nomenclature[H1Sigma]{$H^1(\Sigma^g_{\K} M)$}{first order Sobolev space of sections of $\Sigma^g_{\K} M$}

For this entire article, let $(M, \Theta)$ be a closed Riemannian spin manifold of dimension $m$ and $\Theta:\GLtp M \to \GLp M$ be a fixed topological spin structure on $M$. For any Riemannian metric $g$ on $M$, we denote by $\Sigma^g_{\K} M \to M$ the spinor bundle with respect to $g$ and $\K \in \{\R, \C\}$. The associated Dirac operator is denoted by $\Dirac^g_{\K}$. We think of this operator as an unbounded operator
\begin{align*}
	\Dirac^g_{\K}: H^1(\Sigma^g_{\K} M) \subset L^2(\Sigma^g_{\K} M) \to L^2(\Sigma^g_{\K} M)
\end{align*} 

densely defined on the first order Sobolev space $H^1(\Sigma^g_{\K} M)$ of sections of $\Sigma^g_{\K} M$. In that sense the operator has a spectrum $\spec \Dirac^g_{\K} \subset \R$. One is usually interested in the case $\K=\C$. In terms of a local orthonormal frame, the Dirac operator is given by $\Dirac^g_{\K} = \sum_{i=1}^{m}{e_i \cdot \nabla_{e_i}^g}$ and its spectrum comprises of those $\lambda \in \R$ for which there exists a non-trivial spinor field $\psi \in \Gamma(\Sigma^g_{\K} M) $ such that
\begin{align}
	\label{EqDirac}
	\Dirac^g_{\K} \psi = \lambda \psi.
\end{align}
The equation \Cref{EqDirac} is called the \emph{Dirac equation} and our main result about it is as follows.

\begin{MainThm}[existence of higher multiplicities]
	\label{MainThmHigher}	
	Let $(M, \Theta)$ be a closed spin manifold of dimension $m \equiv 0,6,7 \mod 8$. There exists a Riemannian metric $\tilde g$ on $M$ such that the complex Dirac operator $\Dirac^{\tilde g}_{\C}$ has at least one eigenvalue of multiplicity at least two. In addition, $\tilde g$ can be chosen such that it agrees with an arbitrary metric $g$ outside an arbitrarily small open subset on the manifold.
\end{MainThm}

\subsection{Dahl's conjecture}

The result of \Cref{MainThmHigher} fits nicely into the context of a conjecture by Dahl in \cite{DahlPresc}, which deals with the question of what sequences of real numbers can occur as Dirac spectra. In general, the Dirac spectrum depends on the metric and even on the spin structure, see \cite{FriedTori}. On the other hand, all Dirac spectra have certain properties in common.  

\begin{Lem}[Properties of Dirac spectra]
	\label{LemDiracSpecProp}
	Let $(M,\Theta)$ be a closed spin manifold and $g$ be any Riemannian metric on $M$. Then $\Dirac^g_{\C}$ is a self-adjoint elliptic first order differential operator and its spectrum satisfies the following properties:
	\begin{enumerate}[label=(D\arabic*)]
		\item
			\label{ItDiracPropDiscrete}
			$\spec \Dirac^g_{\K} \subset \R$ is discrete and unbounded from both sides.
		\item
			\label{ItDiracPropQuarternionic}
			In case $m \equiv 2,3,4 \mod 8$, there exists a \emph{quaternionic structure} and all eigenspaces are even-dimensional over $\C$.
			
		\item
			\label{ItDiracPropSymmetric}
			In case $m \not \equiv 3 \mod 4$, the Dirac spectrum is symmetric about zero including multiplicities.

		\item
			\label{ItDiracPropAtiyah}
			The kernel of the Dirac operator satisfies the estimate
			\begin{align*}
				\dim_{\C} \ker \Dirac^g_{\C} \geq  
				\begin{cases}
					|\hat A(M)|, & m \equiv 0,4 \mod 8, \\
					1, & m \equiv 1 \mod 8 \text{ and } \alpha(M) \neq 0, \\
					2, & m \equiv 2 \mod 8 \text{ and } \alpha(M) \neq 0.
				\end{cases}
			\end{align*}
			Here, $\hat A(M)$ denotes the $\hat A$-genus and $\alpha(M)$ denotes the $\alpha$-genus. 
		\item
			\label{ItDiracPropWeyl}
			The growth of the Dirac eigenvalues satisfies a certain Weyl's law. 
	\end{enumerate}
\end{Lem}

For a proof of these elementary facts as well as for an introduction into spin geometry in general, the reader is referred to \cite{LM,FriedSpinGeom,hijazi}. 

\Cref{LemDiracSpecProp} raises the question whether or not one can prescribe Dirac spectra artibrarily as long as one does not violate its assertions.

\begin{Cnj}[\protect{\cite{DahlPresc}}]
	\label{CnjDahl}
	Let $k \in \N$, $\Lambda_1, \Lambda_2 \in \R$, $\Lambda_1 < \Lambda_2$, and $(M,\Theta)$ be a compact spin manifold. For any non-zero $\lambda_1 \leq \ldots \leq \lambda_k \in ]\Lambda_1, \Lambda_2[$ satisfying \cref{ItDiracPropQuarternionic,ItDiracPropSymmetric}, there exists a metric $g$ on $M$ such that
	\begin{align*}
		\spec \Dirac^g_{\C} \cap \mathopen{]} \Lambda_1, \Lambda_2 [ = \{\lambda_1 \leq \ldots \leq \lambda_k\},
	\end{align*}
	where the eigenvalues are counted with multiplicities.
\end{Cnj}

In the same article, Dahl also gives a proof of this conjecture in the case where all eigenvalues are \emph{simple}, see \cite[Thm. 1]{DahlPresc}. In this context, an eigenvalue $\lambda$ is called \emph{simple} if its \emph{multiplicitiy}
\begin{align*}
	\mu(\lambda) := \dim_{K} \ker (\Dirac^g_{\C} - \lambda), && 
	K := \begin{cases}
			\H, & m \equiv 2,3,4 \mod 8, \\
			\C, & \text{otherwise}
		\end{cases}
\end{align*}
is equal to $1$. Here, $\H$ denotes the quaternions.
\nomenclature[K]{$K$}{$\C$ or $\H$}
\nomenclature[mu]{$\mu(\lambda)$}{multiplicitiy of the eigenvalue $\lambda$}

\begin{Rem}[multiplicities]
	Denote by $\mu_{\K}(\lambda) := \dim_{\K}\ker (\Dirac^g_{\K} - \lambda)$ the multiplicity of an eigenvalue $\lambda$ over $\K \in \{\R, \C \}$. Then the various notions of multiplicity are related by
		\begin{align*}
		\mu(\lambda) = 
		\begin{cases}
			\mu_{\C}(\lambda) = \mu_{\R}(\lambda), & m \equiv 0,6,7 \mod 8, \\
			\mu_{\C}(\lambda) = \tfrac{1}{2} \mu_{\R}(\lambda), & m \equiv 1,5 \mod 8, \\
			\tfrac{1}{2} \mu_{\C}(\lambda) = \tfrac{1}{4} \mu_{\R}(\lambda), & m \equiv 2, 3, 4 \mod 8.
		\end{cases}
	\end{align*}
	We will be primarily concerned with the case $m \equiv 0,6,7 \mod 8$, where all these notions agree, see also \Cref{RemRealvsComplex}.
\end{Rem}
\nomenclature[muK]{$\mu_{\K}(\lambda)$}{multiplicitiy of the eigenvalue $\lambda$ over $\K$}

The question what one can say about \Cref{CnjDahl} in case of higher multiplicities has been open ever since. One would guess that one can prescribe eigenvalues of arbitrary finite multiplicity. Unfortunately, the proof of \Cref{CnjDahl} in case of simple multiplicities does not carry over to higher multiplicities. Therefore, the aim of this article is to introduce some new techniques to approach \Cref{CnjDahl} in case of higher multiplicities, which will allow us to prove \Cref{MainThmHigher}.

\begin{Rem}[real vs. complex spin geometry]
	\label{RemRealvsComplex}
	The restriction in the dimension in the assertion of \Cref{MainThmHigher} stems from the fact that we need tools from real and from complex spin geometry. In dimensions $m \equiv 0,6,7 \mod 8$, complex spin geometry is the complexification of real spin geometry. More precisely, the complexification of an irreducible real representation of the real Clifford algebra will be an irreducible complex representation of the complex Cilfford algebra. This follows from the explicit classification of real and complex Clifford algebras, see for instance \cite[I.\textsection 4]{LM}. Hence the complexification of the real spinor representation is a complex one. This behavior under complexification goes through for all other structures on the spinor bundle, in particular Clifford multiplication, the spinorial connection and the Dirac operator. Thus, in dimensions $m \equiv 0,6,7 \mod 8$, we can jump back and forth between the real and the complex spin geometry. 
\end{Rem}

\begin{Rem}[neighborhood]
	The precise nature of the neighborhood mentioned in \Cref{MainThmHigher} will become clear in the proof. It will be a surgery disc around a point, where we perform a connected sum, see in \Cref{FigSphereSurgery}. However, $\tilde g$ will typically not be in a small $\mathcal{C}^1$-neighborhood of $g$ in the space $\Rm(M)$ of Riemannian metrics on $M$.
\end{Rem}

\subsection{Proof strategy}
\label{SubSectPrfStrategyMainThmHigher}

\begin{figure}[t] 
	\begin{center}
		\begin{tikzpicture}

  \tikzset{
    decoration={
      markings,
      mark=between positions 0.125 and 0.90 step 0.0625 with {
        \xdef\maxseq{\pgfkeysvalueof{/pgf/decoration/mark info/sequence number}}
        \coordinate (pt-\maxseq)
        at (0,0);
      },
    }
  }

	
	\draw[black,very thick]
		(0,0) ellipse (4 and 3);

	\draw[gray!10,fill]
		(0.5,0) ellipse (1.5 and 1);

	\draw[black]
		(0.5,0) ellipse (1.5 and 1);

	\draw[thick,postaction={decorate}] (-2.5,-1) to (-2,-0.5)
		to [out=45,in=180] (0,-1.5)
		to [out=0,in=180] (2,-1)
		to [out=0,in=-90] (2.5,0)
		to [out=90,in=0] (1.5,1.5)
		to [out=170,in=0] (0,1.5)
		to [out=180,in=45] (-2,-0.2)
		to [out=220,in=45] (-2.5,-1);
		
	\foreach \c in {1,...,\maxseq}{
	\draw[black,rotate=1.35*360*(\c-1) / \maxseq] (pt-\c) +(0,-0.4) -- +(0,0.4);
	}
	
	\fill[black] (-2.4,-0.9) circle (0.07);
	

	\coordinate[label=above:$Y$] (A) at (-2,1.5);
	\coordinate[label=above:$X \setminus Y$] (A) at (0.5,-0.2);
	\coordinate[label=above:$\gamma$] (gamma) at (-1.7,-1.5);
	
\end{tikzpicture}
		\caption[The ``Lasso Lemma''.]{The ``Lasso Lemma''.}
		\label{FigLasso}
	\end{center}
\end{figure}
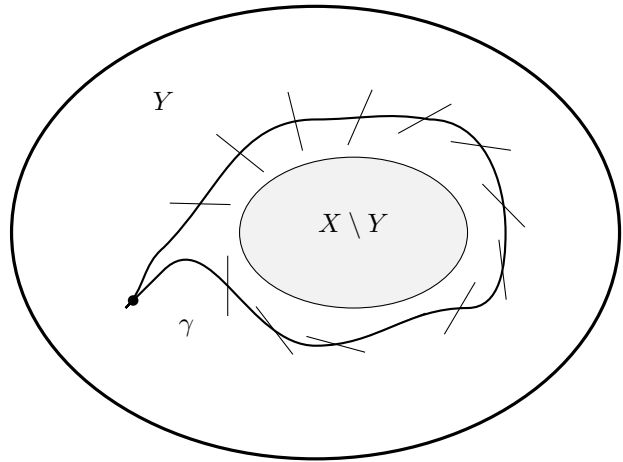

The key idea to prove \Cref{MainThmHigher} is the following simple topological reasoning to which we will refer to as the \emph{Lasso Lemma}, see \Cref{FigLasso}. 

\begin{Lem}[``Lasso Lemma'']
	\label{LemLasso}
	Let $X$ be a simply connected topological space and let $Y \subset X$ be any subspace. Let $\gamma:\S^1 \to Y$ be a loop and $E \to Y$ be a vector bundle such that $\gamma^* E \to \S^1$ is not trivial. Then $X \setminus Y$ is not empty.
\end{Lem}

\begin{Prf}
	Since $X$ is simply connected, there exists a homotopy $H:I^2 \to X$ from $\gamma$ to the constant loop. Since $\gamma^*E \to \S^1$ is not trivial, $\gamma$ cannot be null-homotopic in $Y$. Thus, there has to be at least one point in $X \setminus Y$ that is hit by $H$, hence $X \setminus Y \neq \emptyset$.
\end{Prf}

\begin{Rem}
	\label{RemDiscExtension}
	Of course we can identify $I^2$ with $\D^2$ and obtain that any extension $H:\D^2 \to X$ of $\gamma:\S^1 \to Y$ satisfies $H(p) \in X \setminus Y $ for at least one $p \in \D^2$.
\end{Rem}

We will apply this reasoning in the following way: We set $X(M) := (\Rm(M), \mathcal{C}^1)$, the space of all Riemannian metrics on $M$ endowed with $\mathcal{C}^1$-topology. The set $Y(M)$ will be a subspace of metrics tailor-made such that $X(M) \setminus Y(M) \neq \emptyset$ directly implies the existence of an eigenvalue of higher multiplicity. (The set $Y(M)$ contains the set of all metrics for which all eigenvalues are simple, see \Cref{DefPartEBundleGlobal}. We use $Y(M)$ instead of this simpler set for technical reasons.) The bundle $E := E(M)$ consists of the span of the eigenspinors corresponding to a certain finite set of eigenvalues, see \Cref{DefPartEBundleGlobal}. For the loop $\gamma$ we will have to construct a suitable loop $\mathbf{g}:\S^1 \to Y(M)$ of Riemannian metrics. 

\nomenclature[RmM]{$\Rm(M)$}{space of Riemannian metrics on $M$ with $\mathcal{C}^1$-topology}
\nomenclature[XM]{$X(M)$}{$X(M) = \Rm(M)$}

Unfortunately, we will not be able to construct this loop directly. Therefore, we will use the following strategy: In \Cref{SubSectLoopsMetricsDiffeos}, we consider loops of spin diffeomorphisms $(f_{\alpha})_{\alpha \in \S^1}$ on $M$ and study loops of metrics induced by setting $g_{\alpha} := (f^{-1}_{\alpha})^* g$, $\alpha \in \S^1$, $g \in \Rm(M)$. We will work out a criterion when this loop induces a non-orientable bundle over $\S^1$ as desired, see \Cref{ThmS1ActionBundleTwisted}. This reduces the problem of finding a loop of metrics to finding a loop of spin diffeomorphisms (which might be even harder in general). In \Cref{SubSectSphere}, we will show that the family of rotations by degree $\alpha$ on the sphere $\S^m$ will suit our purpose, if we start with a metric $g_0$ that is obtained from the round metric by a small perturbation. This will give us the desired loop of metrics on the sphere $\S^m$. 

Finally, we will have to transport the loop of metrics on the sphere $\S^m$ to our original manifold $M$. Any smooth $m$-manifold $M$ is diffeomorphic to $M \sharp \S^m$, where $\sharp$ denotes a connected sum, which is a special type of \emph{surgery}. In \Cref{SubSectSurgeryEigenbundles}, we will review the concept of surgery in the setting of Riemannian spin geometry and ultimately show that the existence of a suitable loop of metrics is stable under certain surgeries, see \Cref{ThmSurgeryStabilityOdd}. Applying this to the connected sum will yield the desired result, see also \Cref{FigSphereSurgery}.

\subsection{Acknowledgements}
I would like to thank Bernd Ammann for interesting discussions. I'm also grateful to Mattias Dahl for explaining parts of his previous work to me. This research was kindly supported by the \emph{Studienstiftung des deutschen Volkes} and the \emph{DFG Graduiertenkolleg GRK 1692 ``Curvature, Cycles and Cohomology''}.

\subsection{Comparison to results for Laplace, Schrödinger and other operators}
\label{SectLaplaceSchroedinger}

One should note that \Cref{CnjDahl} has not only been formulated for the Dirac operator. The Laplace operator on functions and the Schrödinger operator has been studies by Colin de Verdière in \cite{verdiere1,verdiere2,verdiere3}. Some parts of these articles are formulated for more general classes of self-adjoint positive operators (notice however that the Dirac operator is not positive). These results were generalized later by Pierre Jammes in \cite{jammes1,jammes3,jammes2,jammes4} to the case of a Hodge Laplacian acting on $p$-forms and even to the Witten Laplacian. 

It is interesting to note how the research on the problem of prescribing the eigenvalues of the Laplace operator has progressed: Jammes started with simple eigenvalues, advanced to double eigenvalues and finally considered eigenvalues of arbitrary multiplicity. Therefore, we think that a similar approach for the Dirac operator is reasonable. 

A similar problem is given by the Laplace operator $\Delta_{\Omega} = - \sum_i \partial_i^2$ on a domain $\Omega \subset \R^m$ with Dirichlet boundary conditions. The spectrum $\{\lambda_j(\Omega)\}_{j \in \N}$ of $\Delta_{\Omega}$ depends on $\Omega$, but it cannot be prescribed arbitrarily by varying $\Omega$ among all domains of $\R^m$ with a fixed volume. By the theorem of Faber-Krahn, the Ball $B$ of volume $c$ satisfies $\lambda_1(B) = \min\{ \lambda_1(\Omega) \mid \Omega \subseto \R^m, |\Omega| = c \}$. Analogously, by the theorem of Kran-Szegö, the minimum of $\lambda_2(\Omega)$ among all bounded open subsets of $\R^m$ with given volume is achieved by the union of two identical balls. A proof of these results (and many more results in this direction) can be found in \cite{henrot}.

While it is possible to prescribe eigenvalues of higher multiplicity for the Laplace operator, there are other physically motivated operators $L$ for which $Lu = \lambda u$ always implies that $\lambda$ is simple. For instance, consider the \emph{Sturm-Liouville operator} $ L u := - \Big( \frac{d}{dx} \left( p \cdot \frac{d}{dx} \right) + q \Big)u = \lambda u$ on $L^2([a,b])$ subject to the boundary conditions
\begin{align}
	\label{EqSturmLiouBoundary}
	c_a u (a) + d_a u'(a) = 0, &&
	c_b u (b) + d_b u'(b) = 0.
\end{align}
for some fixed constants $c_a, d_a, c_b, d_b \in \R$. Here, $p$ is differentiable and positive and $q$ is continuous. As a domain for $L$ we can choose the closure of the $\mathcal{C}^2$ functions satisfying the boundary conditions \Cref{EqSturmLiouBoundary} under the $L^2$-scalar product. Then $L$ is an elliptic self-adjoint operator of second order depending on the functions $p$ and $q$. However, any eigenvalue $\lambda$ of $L$ is always simple regardless of the choice of $p$ and $q$, see for instance \cite[Thm 4.1]{hartmann}.

\section{Construction of the set}

The construction of a subset $Y(M)$ suitable to apply \Cref{LemLasso} needs a consistent enumeration of the spectrum for all metrics. This is possible by the following result.

\begin{Thm}[\protect{\cite[Main Thm. 2]{evpaper}}]
	\label{ThmEvpaper}
	\nomenclature[lambdaj]{$\lambda_j(g)$}{$j$-th eigenvalue of $\Dirac^g_{\K}$}
	There exists a family of continuous functions $\{\lambda_j:\Rm(M) \to \R \}_{j \in \Z}$ such that for all $g \in \Rm(M)$, the sequence $(\lambda_j(g))_{j \in \Z}$ represents all the eigenvalues of $\Dirac^g_{\K}$ (counted with multiplicities) and is non-decreasing, i.e. all $g \in \Rm(M)$ satisfy $\lambda_j(g) \leq \lambda_k(g)$, if $j \leq k$.
\end{Thm}

We fix one such family for the entire article.

\begin{Def}[construction of $Y(M)$]
	\label{DefYKM}
	We define, 
	\begin{align}
		\label{EqDefRmA}
		\begin{split}
			Y(M) := \{g \in \Rm(M) \mid & \exists 1 \leq j \leq k: \mu(\lambda_j(g)) \in 2 \N +1, \\
			& \lambda_{0}(g) < \lambda_1(g), \lambda_{k}(g) < \lambda_{k+1}(g) \},
		\end{split}
	\end{align}
	where $k \in 2 \N +1 $ is a fixed number (whose precise value will be specified later, see \Cref{RemChoicen}).
	\nomenclature[YM]{$Y(M)$}{a subset of Riemannian metrics}
\end{Def}

\begin{Rem}
	One might wonder, why we define $Y(M)$ in such a complicated manner. For the moment we recall that $X(M) = \Rm(M)$ and convince ourselves that 
	\begin{align*}
		X(M) \setminus Y(M)
		& = \{g \in \Rm(M) \mid  \forall 1 \leq j \leq k: \mu(\lambda_j(g)) \in 2 \N \text{ or } \\
		& \phantom{== \{g \in \Rm(M) \mid} \lambda_{0}(g) = \lambda_1(g) \text{ or } \lambda_{k}(g) = \lambda_{k+1}(g) \} \\
		& \subset \{g \in \Rm(M) \mid \exists 1 \leq j \leq k : \mu(\lambda_j(g))  \geq 2 \}.
	\end{align*}
	Therefore, if we can show that $X(M) \setminus Y(M)$ is not empty, we have shown the existence of an eigenvalue of higher multiplicity. 
\end{Rem}

\section{Construction of the bundle}

\subsection{Definition of $E(M)$ as a set}

The construction of the vector bundle $E(M)$ as a set is straightforward.

\begin{Def}[construction of $E(M)$]
	\label{DefPartEBundleGlobal}
	We define
	\begin{align}
		\label{EqDefEM}
		\begin{split}
		E^g(M) &:=  \Lin \{ \psi \in H^1(\Sigma^g_{\R} M) \mid \exists 1 \leq j \leq k: \Dirac^g_{\R} \psi = \lambda_j(g) \psi \} , \\
		E(M) & := \coprod_{g \in Y(M)}{E^g(M)} \to Y(M),
		\end{split}
	\end{align}
	where the bundle projection simply maps a $\psi \in H^1(\Sigma^g_{\R} M)$ to $g$.
\end{Def} \nomenclature[EM]{$E(M)$}{a finite dimensional bundle over $Y(M)$}

\begin{Rem}
	\label{RemWhyEComplicated}
	Recall that $\{\lambda_j\}_{j \in \Z}$ evaluated at any $g \in \Rm(M)$ is a non-decreasing enumeration of the Dirac spectrum $\spec \Dirac^g_{\K}$ counted with multiplicities. This is why we had to add the conditions $\lambda_0(g) < \lambda_1(g)$ and $\lambda_k(g) < \lambda_{k+1}(g)$ in \Cref{EqDefRmA}; they ensure that the vector spaces defined in \Cref{EqDefEM} have constant dimension, thus $E(M)$ has constant rank. 
\end{Rem}

Notice that $E(M)$ consists of real vector spaces, since they are spanned by real eigenspinors of the real Dirac operator. We want to use $E(M)$ to  make a conclusion about the complex Dirac operator, so we will have to jump between the real and the complex spin geometry as discussed in \Cref{RemRealvsComplex}.

\subsection{Topologization of $E(M)$}

It remains only to topologize $E(M)$ and show that it is a continuous vector bundle. The topology will be the subspace topology of a \emph{universal spinor field bundle}. The continuity claim will follow from standard arguments of functional analysis.

For the definition of a topology on $E(M)$, we need to compare the spinors in spinor bundles formed with respect to two different metrics, let's say $g,h \in \Rm(M)$. The problem is that the the two Dirac operators $\Dirac^g_{\K}$ and $\Dirac^h_{\K}$ cannot be compared directly, because not only the operators depend on the metric, but also their domains. Therefore, the expression $\Dirac^g_{\K} - \Dirac^h_{\K}$ does not make any sense. A solution to this problem is to systematically construct identification isomorphisms (of Hilbert spaces)
\begin{align} \nomenclature[barbetagh]{$\bar \beta_{g,h}$}{identification isomorphism $L^2(\Sigma^g_{\K} M) \to L^2(\Sigma^h_{\K} M)$}
	\label{DefDefbarBetagh}
	\bar \beta_{h,g}:L^2(\Sigma^h_{\K} M) \to L^2(\Sigma^g_{\K} M)	
\end{align}
for any two metrics $g$ and $h$ and use these maps to pull back one Dirac operator to the domain of definition of the other. In the Riemannian case, this program has been carried out in \cite{BourgGaud} by means of a connection, but can also be described using only the Lifting Theorem, see \cite{MaierGen}. There is also an alternative approach using \emph{generalized cylinders} that also works in the Lorentz case, see \cite{BaerGaudMor}. We will apply these results in the following way.

\begin{Thm}[universal spinor field bundle]
	\nomenclature[L2SigmaM]{$L^2(\Sigma_{\K} M)$}{universal spinor field bundle}
	\label{ThmUniSpinorFieldBdle} 
	The \emph{universal spinor field bundle} defined by
		\DefMap{L^2(\Sigma_{\K} M) := \coprod_{g \in \Rm(M)}{L^2(\Sigma^g_{\K} M)}}{\Rm(M)}{\psi \in L^2(\Sigma^g_{\K} M)}{g}
	has a unique topology as a Hilbert bundle such that for any $g \in \Rm(M)$, 
	\DefMap{\bar \beta^g:L^2(\Sigma_{\K} M)}{L^2(\Sigma^g_{\K} M) \times \Rm(M)}{\psi \in L^2(\Sigma^h_{\K} M)}{(\bar \beta_{h,g}(\psi),h),}
	is a global trivialization. Here, $\bar \beta_{h,g}$ is the identification isomorphism \Cref{DefDefbarBetagh}.
\end{Thm}

\begin{Prf}
	We fix a metric $g \in \Rm(M)$ and define the topology on $L^2(\Sigma_{\K} M)$ by simply declaring $\bar \beta^g$ to be a trivialization. To see that this topology is independent of $g$, one has to show that the identification isomorphisms $\bar \beta_{h,g}$ themselves depend $\mathcal{C}^1$-continuously on the metric. This is clear from the construction, but a bit tedious to carry out, see \cite[Chapter 4]{NiknoDiss} for details.
\end{Prf}

\begin{Thm}[continuity of eigenbundles]
	\label{ThmContEigenbundle}
	Let $Y \subset \Rm(M)$ be any subspace and $k \in \N$ such that
	\begin{align*}
		\forall g \in Y: \lambda_{0}(g) < \lambda_1(g), \; \lambda_{k}(g) < \lambda_{k+1}(g).
	\end{align*}
	Then the \emph{eigenbundle}
	\begin{align*}
		E := \coprod_{g \in Y}{\sum_{j=1}^{k}{\ker(\Dirac^{g}_{\K} - \lambda_j(g))}} \to Y
	\end{align*}
	is a continuous vector bundle of rank $k$ over $\K$, when endowed with the subspace topology inherited from the universal spinor field bundle $L^2(\Sigma^g_{\K} M)$ from \Cref{ThmUniSpinorFieldBdle}.
\end{Thm}

\begin{Prf}
	For any $g$ in $Y$, we can find a simple closed curve $c:\S^1 \to \C$ such that $\lambda_1(g), \ldots, \lambda_k(g)$ lie inside the area enclosed by $c$ and the rest of the spectrum lies outside this area. Since the $\lambda_j$'s are continuous, the same holds in a small neighborhood of $g$. We obtain that the expression
	\begin{align}
		\label{EqSpecProj}
		P^g(\psi) := -\frac{1}{2 \pi i} \oint_{c}{(z-\Dirac^g_{\K})^{-1}\psi dz} 
	\end{align}
	depends continuously on $g$. It is shown in \cite[Theorem 6.17, p. 178]{kato} that $P^g$ and $\id - P^g$ define operators with spectrum $\lambda_1(g), \ldots, \lambda_k(g)$ respectively $\{\lambda_j(g), j \neq 1, \ldots, k\}$. Since the Dirac operator is self-adjoint, it follows that \cref{EqSpecProj} is actually the spectral projection onto the sum of eigenspaces spanned by $\lambda_1(g), \ldots, \lambda_k(g)$. As a result the images of the various $P^g$'s assemble to a continuous vector bundle, see \cite[Thm. 4.5.2]{NiknoDiss} for more details.
\end{Prf}

\begin{Cor}[topologization of $E(M)$]
	\label{ThmEContBundle}
	The bundle $E(M) \to Y(M)$ from \Cref{DefPartEBundleGlobal} is a continuous vector bundle of rank $k$ 
\end{Cor}

\subsection{Triviality of vector bundles over $\S^1$}
\label{SubSectTrivVB}

Ultimately, we want to apply \Cref{LemLasso} and therefore, we will have to verify that a real vector bundle over $\S^1$ is not trivial. The question whether or not a vector bundle is trivial can in general be approached by various topological machineries. But we are mainly interested in vector bundles over $\S^1$ and here the situation is very simple: The set of isomorphism classes of vector bundles of rank $k$ over $\S^1$ has only two elements, see for instance \cite[p.25]{hatcher}. One class represents the trivial, hence orientable bundles, the other class consists of vector bundles that are non-orientable, hence non-trivial. 

\begin{Rem}[sign of a vector bundle]
	\label{RemFrameCurves}
	A neat criterion to check when a real vector bundle $E \to \S^1$ of rank $k$ is not orientable is the following: Let $\GL E \to \S^1$ be the principal $\GL_k$-bundle of frames of $E$. Let $I := [0,1]$ be the unit interval and denote by $\pi_{\S^1}:I \to \S^1$ the canonical projection. Since $I$ is contractible, $\pi_{\S^1}^*(\GL E) \to I$ has a global section $\Psi$. For any such section $\Psi$, there exists $A \in \GL_k$ such that $\Psi(1)=\Psi(0).A$. Clearly, 
	\begin{align*}
		\det(A) > 0 && 
		\Longleftrightarrow &&  \text{$E$ is orientable}
		&&\Longleftrightarrow && \text{$E$ is trivial}.
	\end{align*}
	We define $\sgn(E) := \sgn(\Psi) := \sgn(\det(A)) \in \Z_2:=\{\pm 1 \}$ to be the \emph{sign of $E$}.
	\nomenclature[I]{$I$}{$I := [0,1] \subset \R$ is the unit interval}
	\nomenclature[piS1]{$\pi_{\S^1}:I \to \S^1$}{canonical projection}
	\nomenclature[sgnE]{$\sgn(E)$}{sign of a vector bundle}
	\nomenclature[Z2]{$\Z_2$}{$\{\pm 1\}$}
\end{Rem}

It will be very important that the sign of a vector bundle is stable under small deformations of the bundle in the following sense.

\begin{Thm}[sign stability]
	\label{ThmSignStability}
	Let $\mathcal{H} \to X$ be a Hilbert bundle and $E,\tilde E \to X$ be two $k$-dimensional subbundles of $\mathcal{H}$ with induced metric. Denote by $S \tilde E \to X$ the bundle of unit spheres of $\tilde E$. If
	\begin{align}
		\label{EqSignStabDistHyp}
		\forall x \in X: \dist(E_{x}, S \tilde E_{x}) < 1,
	\end{align}
	then $E \cong \tilde E$. In particular, if $X=\S^1$, then $\sgn(E) = \sgn(\tilde E)$.
\end{Thm}

\begin{Prf}
	Let $P:\mathcal{H} \to \mathcal{H}$ be the orthogonal projection onto $E$. Let $x \in X$, $\tilde v \in S \tilde E_x$ be arbitrary and assume $P_x(\tilde v) = 0$. By definition, this simply means that $\tilde v$ is perpendicular to $E_{x}$. This implies
	\begin{align*}
		\dist(E_{x}, \tilde v) 
		=\|P_{x}(\tilde v) - \tilde v\|
		=\|\tilde v\|
		=1,
	\end{align*}
	which contradicts our assumption \Cref{EqSignStabDistHyp}. Consequently, $P|_{\tilde E}: \tilde E \to E$ is an isomorphism. 
\end{Prf}

\section{Construction of the loop}
\label{SectConstrLoop}

\subsection{Loops of metrics via loops of diffeomorphisms}
\label{SubSectLoopsMetricsDiffeos}

In this section, we introduce a technique to produce certain loops of metrics via loops of spin diffeomorphisms. We denote by $\Diff(M)$ the diffeomorphism group of $M$ endowed with the usual $\mathcal{C}^{\infty}$-topology, see for instance \cite[Chpt. 2.1]{HirschDiffTop}. We will also use this topology on all the other mapping spaces. 

\nomenclature[DiffM]{$\Diff(M)$}{the diffeomorphism group of $M$}

\begin{Def}[associated loops of metrics]
	\label{DefAssocSpinLoop}
	Let $f:\S^1 \to \Diff(M)$ be a loop of diffeomorphisms and $g \in \Rm(M)$ be any Riemannian metric. The family of metrics 
		\DefMap{\mathbf{g}:\S^1}{\Rm(M)}{\alpha}{g_\alpha := (f_\alpha^{-1})^*g}
	is called an \emph{associated loop of metrics}. 
\end{Def}

We are primarily interested in loops spin diffeomorphisms. Since there exist slightly different conventions, we fix the following notion. 

\begin{Def}[spin diffeomorphism]
	An orientation-preserving diffeomorphism $f$ of $M$ is a \emph{spin diffeomorphism} (or just ``\emph{is spin}''), if there exists $\hat f$ such that
	\begin{align*}
		\xymatrix{
			\GLtp M
				\ar@{-->}[r]^-{\hat f}
				\ar[d]^-{\Theta}_{2:1}
			& \GLtp M
				\ar[d]^-{\Theta}_{2:1}
			\\
			\GLp M
				\ar[r]^-{f_*}
				\ar[d]
			&\GLp M
				\ar[d]
			\\
			M
				\ar[r]^-{f}
			&M
		}
	\end{align*}
	commutes. We say $\hat f$ is a \emph{spin lift of $f$}. We define
	\begin{align*}
		\Diff^{\spin}(M) & := \{f \in \Diff(M) \mid \text{$f$ is spin }\}, \\
		\widehat{\Diff}^{\spin}(M) &:= \{(f,\hat f) \mid \Theta \circ \hat f = f_* \circ \Theta \},
	\end{align*}
	the \emph{spin diffeomorphism group (with lift)} .
	\nomenclature[DiffspinM]{$\Diff^{\spin}(M)$}{group of spin diffeomorphisms}
	\nomenclature[DiffhatspinM]{$\widehat{\Diff}^{\spin}(M)$}{group of spin diffeomorphisms with lift}
\end{Def}
Notice that in case $M$ is connected and $f$ is spin, there are always two spin lifts $\hat f_{\pm}$ of $f$ and $\hat f_{-} = \hat f_{+}.(-1)$, where $.(-1)$ denotes the action of $-1 \in \GLtp_m$. With a bit more work, one can show the following relation.

\begin{Thm} 
	Let $M$ be connected. The canonical projection
	\begin{align*}
		\pr^{\spin}: \widehat{\Diff}^{\spin}(M) \to \Diff^{\spin}(M), &&
		(f,\hat f) \mapsto f,		
	\end{align*}
	is a $2:1$ covering space.
	\nomenclature[prspin]{$\pr^{\spin}$}{projection $\widehat{\Diff}^{\spin}(M) \to \Diff^{\spin}(M)$}
\end{Thm}

\begin{Prf}
	This follows essentially from the fact that $\Theta:\GLtp M \to \GLp M$ is a $2:1$-covering and that locally $f_* = \Theta \circ \hat f \circ \Theta^{-1}$. A detailed proof can be found in \cite[Thm. 2.6.4]{NiknoDiss}.
\end{Prf}

\begin{Def}[odd/even]
	\label{DefOddEven}
	A loop of spin diffeomorphisms $f:\S^1 \to \Diff^{\spin}(M)$ is \emph{even}, if there exists a loop $\hat f$ such that
	\begin{align}
		\label{EqDefOddEven}
		\begin{split}
			\xymatrix{
				& \widehat{\Diff}^{\spin}(M)
					\ar[d]^{\pr^{\spin}}_-{2:1}
				\\
				\S^1
					\ar[r]^-{f}
					\ar@/^1pc/@{-->}[ur]^-{\hat f}
				& \Diff^{\spin}(M)
			}
		\end{split}
	\end{align}
	commutes. A loop is \emph{odd}, if it is not even.
\end{Def}

\begin{Rem}[associated isotopy]
	\label{RemAssocIsotopy}
	Recall that we think of $\S^1$ as $I / \sim$ and that $\pi_{\S^1}:I \to \S^1$ denotes the canonical projection. Clearly, $h := f \circ \pr^{\spin}:I \to \Diff^{\spin}$ is a path. By the path lifting property of covering spaces, there exists a lift $\hat h:I \to \widehat{\Diff}^{\spin}(M)$ of $h$. Nevertheless, $\hat h$ will in general only be a path, but not a loop:
	\begin{align}
		\label{EqDefhpath}
		\begin{split}
			\xymatrix{
				& \widehat{\Diff}^{\spin}(M)
					\ar[d]^{\pr^{\spin}}_-{2:1}
				\\
				I
					\ar[r]^-{h}
					\ar@{-->}@/^1pc/[ur]^-{\hat h}
				& \Diff^{\spin}(M).
			}
		\end{split}
	\end{align}
	We say $h$ is the \emph{isotopy associated to $f$}.
\end{Rem}

\begin{Def}[sign]
	Let $f:\S^1 \to \Diff^{\spin}(M)$ be a loop and $h$ be its associated isotopy. The unique number $\sgn(f) \in \Z_2$ such that $\hat h(0) = \sgn(f) \hat h(1)$ is called the \emph{sign of $f$}. 
\end{Def}
\nomenclature[sgnf]{$\sgn(f)$}{sign of a loop}

The sign of $f$ does not depend on the choice of the lift $\hat h$. Apparently, $f$ is even if and only if $\sgn(f)=+1$. One can show that the sign has the following abstract characterization.

\begin{Lem}
	\label{LemSgnAbstract}
	The map $\pr^{\spin}:\widehat{\Diff}^{\spin} \to \Diff^{\spin}(M)$ is a principal $\Z_2$-bundle. The connecting homomorphism $\delta$ from its long exact homotopy sequence
	\begin{align*}
		\xymatrix{
			\pi_1(\widehat{\Diff}^{\spin}(M), (\id_{M}, \id_{\GLtp M}))
			\ar[r]^-{\pr^{\spin}_{\sharp}}
			& \pi_1(\Diff^{\spin}(M), \id_{M})
			\ar[r]^-{\delta}
			& \pi_0(\{ ( \id_{M}, \pm {\id}_{\GLtp M} ) \}) 
		}
	\end{align*}
	satisfies $\delta(f) = (\id_M, \sgn(f) \id_{\GLtp M})$ for any loop $f:(\S^1,0) \to (\Diff^{\spin}(M), \id_M)$.
\end{Lem}

\begin{Prf}
	Any $2:1$-covering is normal, hence a principal $\Z_2$-bundle. Therefore, the claim follows from the definition of the connecting homomorphism $\delta$.
\end{Prf}

\begin{Lem}
	\label{LemSgnGrpHom}
	The sign induces a group homomorphism
	\begin{align*}
		\sgn:\pi_1(\Diff^{\spin}(M), \id_M) \to \Z_2
	\end{align*}
	and $\ker \sgn$ are precisely the homotopy classes of even loops.
\end{Lem}

\begin{Prf}
	It follows from \Cref{LemSgnAbstract} that $\sgn$ is well-defined on homotopy classes. To see that $\sgn$ is a group homomorphism, let $f^{(1)}, f^{(2)} \in \pi_1(\Diff^{\spin}(M), \id_M)$ and consider
		\DefMap{f := f^{(2)} * f^{(1)}:\S^1}{\Diff^{\spin}(M)}{t}{ 
			\begin{cases}
				f^{(1)}(2t), & 0 \leq t \leq \tfrac{1}{2}, \\
				f^{(2)}(2t-1), & \tfrac{1}{2} \leq t \leq 1.
			\end{cases}
		}
	Let $\hat f^{(1)}$ be a lift of $f^{(1)}$ starting at the identity and $\hat f^{(2)}$ be a lift of $f^{(2)}$ starting hat $\hat f^{(1)}(1)$. Then $\hat f = \hat f^{(2)} * \hat f^{(1)}$ is a lift of $f$ and 
	\begin{align*}
		\hat f(1)
		&=\hat f^{(2)}(1)
		=\sgn(f^{(2)}) \hat f^{(2)}(0) 
		=\sgn(f^{(2)}) \hat f^{(1)}(1) \\
		&=\sgn(f^{(2)}) \sgn(f^{(1)}) \hat f^{(1)}(0)
		=\sgn(f^{(2)}) \sgn(f^{(1)}) \hat f(0),
	\end{align*}
	thus $\sgn(f) = \sgn(f^{(2)}) \sgn(f^{(1)})$. Clearly, the constant map $\S^1 \to \Diff^{\spin}(M)$, $\alpha \mapsto \id_M$, lifts to the constant map $\alpha \mapsto (\id_M, \id_{\GLtp M})$, so $\sgn$ is a group homomorphism as claimed.
\end{Prf}

\begin{Rem}
	Although one cannot just replace the isotopy $h$ associated to the loop $f$ by the loop $f$ itself in \Cref{EqDefOddEven}, there always exists a lift $\tilde f$ such that 
	\begin{align*}
		\begin{split}
			\xymatrix{
				\S^1
					\ar@{-->}[r]^-{\tilde f}
					\ar[d]^-{\cdot 2}
				& \widehat{\Diff}^{\spin}(M)
					\ar[d]^-{\pr^{\spin}}_{2:1}
				\\
				\S^1
					\ar[r]_-{f}
				& \Diff^{\spin}(M)
			}
		\end{split}
	\end{align*}
	commutes. Here, $\cdot 2$ denotes the non-trivial double cover of $\S^1$. This follows from the fact that due to \Cref{LemSgnGrpHom}, we have $\sgn(f \circ \cdot 2) = \sgn(f)^2 = +1$ and thus, $f \circ \cdot 2$ is even (although $f$ itself might not be even).
\end{Rem}

\begin{Rem}
	In the special case where $f:\S^1 \to \Diff^{\spin}(M)$ is a group action, the notions of odd and even in the sense of \Cref{DefOddEven} coincide with the notions of an odd and even group action in the sense of \cite[IV.\textsection 3, p.295]{LM}.
\end{Rem}

\begin{Rem}
	\label{RemAssocMetricsSpinIsom}
	Let $f:\S^1 \to \Diff^{\spin}(M)$ be a loop and $h$ be its associated isotopy as in \Cref{RemAssocIsotopy}. For any $t \in I$, we get an induced isometry
		\DefMap{\bar h_t: \Sigma^{g_0}_{\K} M}{\Sigma^{g_t}_{\K} M}{{\psi=[s,v]}}{{[\hat h_t(s),v]}}
	between all the spinor bundles $\Sigma^{g_t}_{\K} M$. The induced map on sections, denoted by $\bar h_{t}$, satisfies
	\begin{align}
		\label{EqSpinDiffeoDiracCommute}
		\Dirac^{g_{t}}_{\K} \circ \bar h_{t} = \bar h_{t} \circ \Dirac^{g_0}_{\K}
	\end{align}
	and therefore maps eigenspinors to eigenspinors.
\end{Rem}

The following will be crucial to verify the hypothesis of \Cref{LemLasso}.

\begin{Thm}
	\label{ThmS1ActionBundleTwisted}
	Let $Y \subset \Rm(M)$ be any subset, $f:\S^1 \to \Diff^{\spin}(M)$ be a loop of spin diffeomorphisms and $g \in \Rm(M)$ such that the associated loop of metrics $\alpha \mapsto (f_{\alpha}^{-1})^* g$ is a map $\mathbf{g}:\S^1 \to Y$. Furthermore, let $E \subset L^2(\Sigma_{\R} M) \to Y$ be a vector bundle of rank $k$. Let $\bar h_t$ be the map induced by $f$ as in \Cref{RemAssocMetricsSpinIsom} and assume $\bar h_t(E) \subset E$ for any $t \in I$. Then $\mathbf{g}^*E \to \S^1$ is not orientable if and only if $f$ is odd and $k$ is odd, i.e. 
	\begin{align*}
		\sgn(\mathbf{g}^*E) =
		\begin{cases}
			-1, & \text{$f$ is odd and $k$ is odd}, \\
			+1, & \text{otherwise},
		\end{cases}
	\end{align*}
	where $\sgn$ is as in \Cref{RemFrameCurves}.
\end{Thm}

\begin{Prf} 
	For any basis $(0,(\psi_1, \ldots, \psi_k)) \in \mathbf{g}^*E|_0$, the curve 
		\DefMap{\Psi:I}{\GL(\mathbf{g}^*E)}{t}{(t,(\bar h_{t}(\psi_1), \ldots, \bar h_{t}(\psi_n)))}
	is a curve of frames for $\mathbf{g}^*E \to \S^1$ as in \Cref{RemFrameCurves}. By definition, we have $\bar h_1 = \sgn(f) \bar h_0$. Consequently, $\Psi(1)=\Psi(0).A$, where $A = \sgn(f) \I_k$, which has determinant $\sgn(f)^{k}$. By definition,
	\begin{align*}
		\sgn(\mathbf{g}^*E)
		=\sgn(\Psi) 
		=\sgn(\det(A))
		= \sgn(f)^k,
	\end{align*}
	which implies the result.
\end{Prf}

\subsection{The sphere}
\label{SubSectSphere}
We have not yet shown that there exists an odd loop $f:\S^1 \to \Diff^{\spin}(M)$ and in general it is very difficult to construct non-trivial loops of spin diffeomorphisms. Fortunately, the most obvious candidate on the sphere does the job. 

\begin{Thm}
	\label{ThmRotationsOddLoop}
	For any $\alpha \in \R$, $m \in \N$, we define the rotation
	\begin{align*}
		R_{\alpha}:=
		\begin{pmatrix}
			\I_{m-1} & 0 & 0 \\
			0 & \cos(\alpha) & -\sin(\alpha) \\
			0 & \sin(\alpha) & \cos(\alpha) 
		\end{pmatrix}:\R^{m+1} \to \R^{m+1}.
	\end{align*}
	The map $f:\S^1 \to \Diff^{\spin}(\S^m)$, $\alpha \mapsto R_{2 \pi \alpha}|_{\S^m}$,  is an odd loop of spin diffeomorphisms.	
	\nomenclature[Ralpha]{$R_{\alpha}$}{rotation by an angle $\alpha$}
\end{Thm}

\begin{Prf}
	Chose the round metric $g\degree$ on $\S^m$. In this case, the spin structure on $\S^m$ is simply given by the universal cover $\vartheta_{m+1}: \Spin_{m+1} \to \SO_{m+1}$. By a tedious calculation carried out in \cite[Lem. 5.4.5]{NiknoDiss}, one can check that the lift $\hat f$ of $f$ is given by
	\begin{align*}
		\forall \alpha \in \S^1: v \mapsto (\cos(\tfrac{\alpha}{2}) + \sin(\tfrac{\alpha}{2})e_{m-1}e_m)v.
	\end{align*}
	It follows from this explicit formula that $f$ is odd.
\end{Prf}

To obtain an associated loop of metrics $g_{\alpha} = (f^{-1}_{\alpha})^*g_0$, we need a start metric $g_0$. Obviously, we cannot take the round metric $g\degree$, since rotations are an isometry with respect to $g\degree$, so the resulting loop would be trivial. A way out is provided by the following. 

\begin{Thm}[odd neighborhood theorem]
	\label{ThmNbhdRoundSphereOdd}
	Let $(M,\Theta)$ be a closed spin manifold of dimension $m \equiv 0,6,7 \mod 8$ and $g_0$ be any Riemannian metric on $M$. In every $\mathcal{C}^1$-neighborhood of $g_0 \in \Rm(M)$, there exists $g \in \Rm(M)$ such that $\Dirac^g_{\K}$ has an eigenvalue $\lambda$ of odd multiplicity. 
\end{Thm}
		
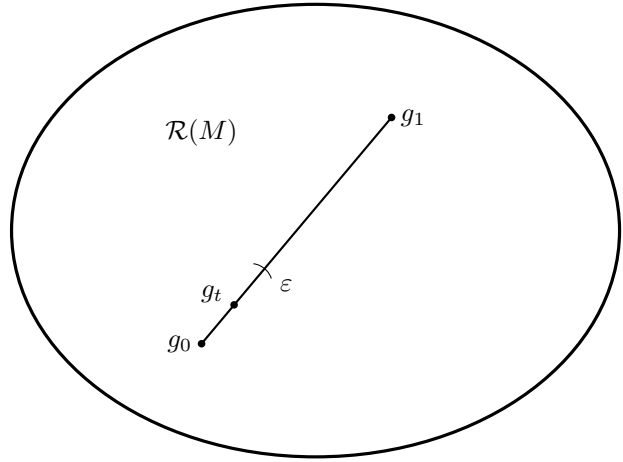
\begin{figure}[t] 
	\begin{center}
		\begin{tikzpicture}

\tikzset{
	partial ellipse/.style args={#1:#2:#3}{
		insert path={+ (#1:#3) arc (#1:#2:#3)}
	}
}

\tikzset{
	decoration={
		markings,
			mark=between positions 0.3 and 0.90 step 1 with {
			\draw[thin] (0,0) [partial ellipse=50:-50:0.125 and 0.2];
			\coordinate[label=right:$\varepsilon$] (eps) at (-0,-0.2);
			\fill[black] (-0.5,0) circle (0.05);
			\coordinate[label=left:$g_t$] (t) at (-0.4,0.1);
		},
	}
}

	\draw[black,very thick]
		(0,0) ellipse (4 and 3);

	\draw[thick,postaction={decorate}] (-1.5,-1.5) to (1,1.5);

	\fill[black] (-1.5,-1.5) circle (0.05);
	\fill[black] (1,1.5) circle (0.05);

	\coordinate[label=above:$\mathcal{R}(M)$] (RM) at (-1.5,1);
	\coordinate[label=left:$g_0$] (A) at (-1.5,-1.5);
	\coordinate[label=right:$g_1$] (gamma) at (1,1.5);

\end{tikzpicture}
		\caption{Finding an odd metric near $g_0$.}
		\label{FigOddnbhd}
	\end{center}
\end{figure}

\begin{Prf}
	The idea of this proof is as follows: By \cite[Thm. 1]{DahlPresc}, there exists a metric $g_1 \in \Rm(M)$ such that $\Dirac^{g_1}_{\K}$ has an eigenvalue of multiplicity $1$, which is odd. Connect the metric $g_0$ with $g_1$, i.e. define the path $g_t := tg_1 + (1-t)g_0$, $t \in I$, see \Cref{FigOddnbhd}. This path is real-analytic. As explained in \cite[Lem. A.0.16]{AndreasDiss}, the Dirac operators $\Dirac^{g_0, g_t}$ are the restriction of a self-adjoint holomorphic family of type (A) onto $I$. Therefore, the eigenvalues of $\Dirac^{g_t}_{\K}$ can be described by a real-analytic family of functions $\{ \lambda_j:I \to \R \}_{j \in \N}$. This means that for any $t \in I$, the sequence $(\lambda_j(t))_{j \in \N}$ represents all the eigenvalues of $\Dirac^{g_t}_{\K}$ counted with multiplicity (but possibly not ordered by magnitude).
	To prove the claim, we argue by contradiction: If the claim is wrong, there exists an open neighborhood around $0$ in which all metrics $g_t$ have only eigenvalues of even multiplicity. Since the eigenvalue functions $\lambda_j$ are real analytic, this behavior extends to all of $I$ and therefore, all the eigenvalues of $g_1$ have even multiplicity as well. But this contradicts the choice of $g_1$. The technical details of this last argument are a bit cumbersome and can be found in \cite[Thm. 5.4.6]{NiknoDiss}.
\end{Prf}

\begin{Rem}[choice of $k$]
	\label{RemChoicen}
	We define the number $k$ to be the dimension of the eigenspace of the eigenvalue of odd multiplicitiy, whose existence is asserted by \Cref{ThmNbhdRoundSphereOdd}.
\end{Rem}

\begin{Rem}
	By \Cref{ThmNbhdRoundSphereOdd}, \Cref{ThmS1ActionBundleTwisted} and 
	\Cref{ThmS1ActionBundleTwisted}, we have verified the hypothesis of \Cref{LemLasso} on the sphere, i.e. we have verified a sufficient criterion, which implies that the sphere admits a metric for which at least one eigenvalue is of higher multiplicity (which is well known). We will show in the next section that this criterion is stable under certain \emph{surgeries}, which will allow us to verify it on much more general manifolds than just the sphere.
\end{Rem}

\subsection{Surgery and eigenbundles}
\label{SubSectSurgeryEigenbundles}
We introduce some basic notions concerning the surgery theory of spin manifolds and recall some well known results by Bär and Dahl published in \cite{BaerDahlSurgery}. Similar techniques are also used in \cite{AmmannSurgery}. We denote by $\S^l$ the $l$-dimensional unit sphere and by $\D^l$ the open unit ball.

\nomenclature[D]{$\D^l$}{Euclidean unit disc}
\nomenclature[S]{$\S^l$}{Euclidean unit sphere}

\begin{DefI}[surgery]
	\label{DefNormalSurgery}
	Let $N$ be a smooth $n$-manifold, let $f:\S^l \times \overline{\D^{n-l}} \to N$, $0 \leq l \leq n$, be a smooth embedding and set $S:=f(\S^l \times \{0\})$, $U := f(\S^l \times \D^{n-l})$. The manifold 
	\begin{align*}
		\tilde N := \Big( (N \setminus U) \amalg (\overline{\D^{l+1}} \times \S^{n-l-1})  \Big)/ \sim,
	\end{align*}
	where $\sim$ is the equivalence relation generated by 
	\begin{align*}
		\forall x \in \S^l \times S^{n-l-1}: x \sim f(x) \in \partial U,
	\end{align*}
	\emph{is obtained by surgery in dimension $l$ along $S$ from $N$}. The number $n-l$ is the \emph{codimension} of the surgery. The map $f$ is the \emph{surgery map} and $S$ is the \emph{surgery sphere}.
\end{DefI}

\begin{figure}[H] 
	\begin{center}
		\begin{tikzpicture}[
  decoration={
    markings,
    mark=between positions 0 and 1 step 0.125 with {\draw[double] (0,-0.1)--(0,0.1); },
  }
]

	\draw[black,thick]
		(0,0) ellipse (4 and 2);

	\draw
		(0,0) ellipse (1 and 0.5);
	\draw[black,thick,postaction={decorate}]
		(0,0) ellipse (1 and 0.5);

	\draw[black,thick,postaction={decorate}]
		(0,3) ellipse (1 and 0.5);		
		
	\draw (-1,3) to (-1,4)
	to [out=90, in=-90] (-2,5)
	to [out=90,in=180] (0,6.5)
	to [out=0,in=90] (2,5)
	to [out=-90, in=90] (1,4)
	to (1,3);

	\coordinate[label=above:$\tilde U$] (Ut) at (0,5);
	\coordinate[label=above:$N \setminus U$] (V) at (-2.5,-0.25);
	\coordinate[label=above:$U$] (U) at (0,-0.25);
	\coordinate[label=above:$S$] (S) at (1,0.25);

\end{tikzpicture}
		\caption[Manifold after surgery.]{The manifold after surgery. Notice that $\partial \tilde U$ is identified with $\partial U \subset (N \setminus U)$.}
		\label{FigSurgeryDec}
	\end{center}
\end{figure}
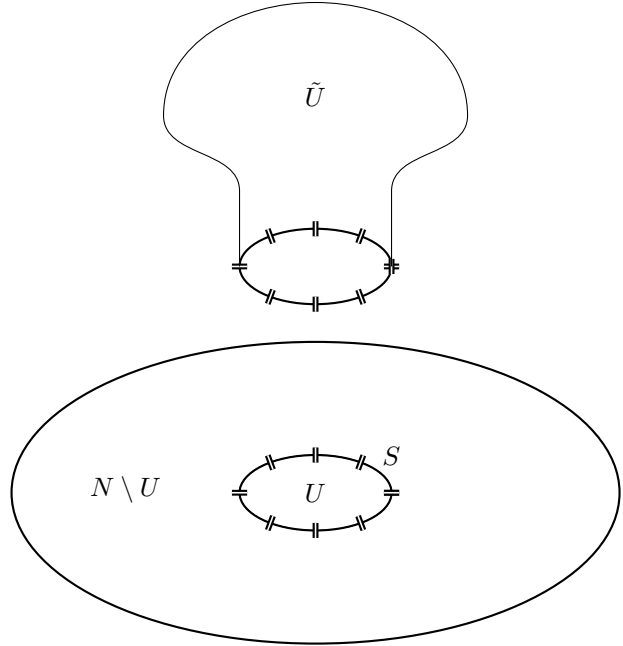

\begin{Rem}
	The space $\tilde N$ is again a smooth manifold (see for instance \cite[IV.1]{kosinski} for a very detailled discussion of the connected sum). The manifold $\tilde N$ is always of the form
	\begin{align}
		\label{EqSurgeryDecompDahl}
		\tilde N = \left( N \setminus U \right) \dot \cup \tilde U,
	\end{align}
	where $\tilde U \subset \tilde N$ is open. Here, by slight abuse of notation, $(N \setminus U) \subset N$ also denotes the image of $N \setminus U$ in the quotient $\tilde N $, see \Cref{FigSurgeryDec}.
\end{Rem}

\begin{Rem}[spin structures and surgery]
	It can be shown that if one performs surgery in codimension $n-l \geq 3$, the spin structure on $N$ always extends uniquely (up to equivalence) to a spin structure on $\tilde N$, if $l \neq 1$. In case $l=1$, the boundary $S^l \times S^{n-l-1}$ has two different spin structures, but only one of them extends to $\D^{l+1} \times \S^{n-l-1}$. Adopting the convention from \cite[p. 56]{BaerDahlSurgery}, we assume that the map $f$ is chosen such that it induces the spin structure that extends. Also, we will only perform surgeries in codimension $n-l \geq 3$.
\end{Rem}

It is natural to ask how the Dirac spectra of a spin manifold before and after surgery are related to one another. The result is roughly that a finite part of the spectrum before surgery is arbitrarily close to the spectrum after surgery. For a precise statement, the following notion is useful.

\begin{Def}[$(\Lambda_1, \Lambda_2, \varepsilon)$-spectral close]
	\index{spectral close}
	Let $T: H \to H$ and $T':H' \to H'$ be two densely defined operators on Hilbert spaces $H$ and $H'$ (over $\K$) with discrete spectrum. Let $\varepsilon > 0$ and $\Lambda_1, \Lambda_2 \in \R$, $\Lambda_1 < \Lambda_2$. Then $T$ and $T'$ are \emph{\emph{$(\Lambda_1, \Lambda_2, \varepsilon)$}}-spectral close if
	\begin{enumerate}
		\item 
			$ \Lambda_1, \Lambda_2 \notin (\spec T \cup \spec T')$.
		\item
			The operators $T$ and $T'$ have the same number $k$ of eigenvalues in $\mathopen{]} \Lambda_1, \Lambda_2 \mathclose{[}$, counted with $\K$-multiplicities.
		\item
			If $\{\lambda_1 \leq \ldots \leq \lambda_k \}$ are the eigenvalues of $T$ in $\mathopen{]} \Lambda_1, \Lambda_2 \mathclose{[}$ and $\{\lambda_1' \leq \ldots \leq \lambda_k' \}$ are the eigenvalues of $T'$ in $\mathopen{]} \Lambda_1, \Lambda_2 \mathclose{[}$, then
			\begin{align*}
				\forall 1 \leq j \leq k: |\lambda_j - \lambda_j'| < \varepsilon.
			\end{align*}
	\end{enumerate}
\end{Def}

Using this terminology, a central result is the following

\begin{Thm}[\protect{\cite[Thm. 1.2]{BaerDahlSurgery}}]
	\label{ThmBaerDahlSurgeryPimped}
	Let $(N^n, g, \Theta^g)$ be a closed Riemannian spin manifold, let $0 \leq l \leq n-3$ and $f:\S^l \times \overline{\D^{n-l}} \to N$ be any surgery map with surgery sphere $S$ as in \Cref{DefNormalSurgery}. For any $\varepsilon > 0$ (sufficiently small) and any $\Lambda > 0$, $\pm \Lambda \notin \spec \Dirac^g_{\C}$, there exists a Riemannian spin manifold $(\tilde N^{\varepsilon}, \tilde g^{\varepsilon})$, which is obtained from $(N, g)$ by surgery such that $\Dirac^g_{\C}$ and $\Dirac^{\tilde g^{\varepsilon}}_{\C}$ are $(-\Lambda, \Lambda, \varepsilon)$-spectral close. This manifold is of the form $\tilde N^{\varepsilon} = (N \setminus U_{\varepsilon}) \dot \cup \tilde U_{\varepsilon}$, where $U_{\varepsilon}$ is an (arbitrarily small) neighborhood of $S$ and the metric $\tilde g^{\varepsilon}$ can be chosen such that $\tilde g|_{N \setminus U_{\varepsilon}} = g|_{N \setminus U_{\varepsilon}}$.
\end{Thm}

We will need not only the statement of \Cref{ThmBaerDahlSurgeryPimped}, but also some arguments from the proof, which relies on estimates of certain \emph{Rayleigh quotients} and these are very useful in their own right. One of the technical obstacles here is that the spinors on $N$ and the spinors on $\tilde N^{\varepsilon}$ cannot be compared directly, since they are defined on different manifolds. A simple yet effective tool to solve this problem are cut-off functions adapted to the surgery.

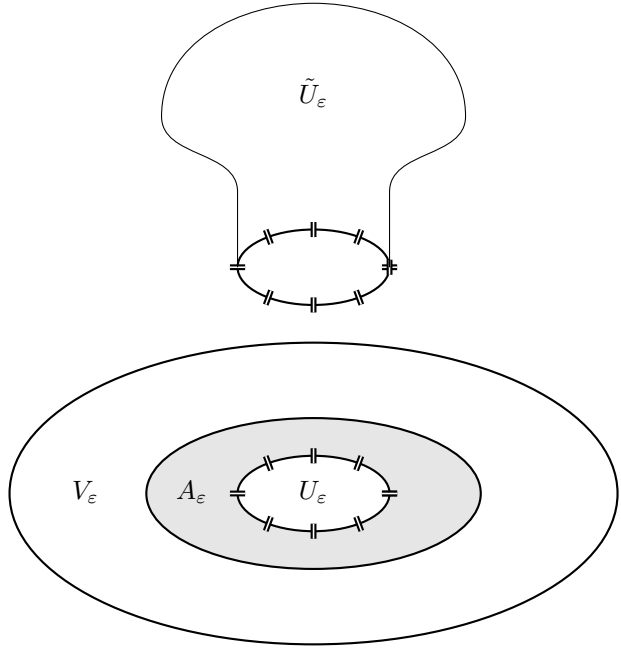
\begin{figure}[t] 
	\begin{center}
		\begin{tikzpicture}[
  decoration={
    markings,
    mark=between positions 0 and 1 step 0.125 with {\draw[double] (0,-0.1)--(0,0.1); },
  }
]

	\draw[black,thick]
		(0,0) ellipse (4 and 2);
		
	\draw[fill,gray!20]
		(0,0) ellipse (2.2 and 1);
	\draw[thick,black]
		(0,0) ellipse (2.2 and 1);

	\draw[fill,white]
		(0,0) ellipse (1 and 0.5);
	\draw[black,thick,postaction={decorate}]
		(0,0) ellipse (1 and 0.5);

	\draw[black,thick,postaction={decorate}]
		(0,3) ellipse (1 and 0.5);		
		
	\draw (-1,3) to (-1,4)
	to [out=90, in=-90] (-2,5)
	to [out=90,in=180] (0,6.5)
	to [out=0,in=90] (2,5)
	to [out=-90, in=90] (1,4)
	to (1,3);

	\coordinate[label=above:$\tilde U_{\varepsilon}$] (Ut) at (0,5);
	\coordinate[label=above:$V_{\varepsilon}$] (V) at (-3,-0.25);
	\coordinate[label=above:$A_{\varepsilon}$] (A) at (-1.6,-0.25);
	\coordinate[label=above:$U_{\varepsilon}$] (U) at (0,-0.25);

\end{tikzpicture}
		\caption[Preparing a manifold for surgery.]{Preparing a manifold $N = U_{\varepsilon} \cup A_{\varepsilon} \cup V_{\varepsilon}$ for surgery.}
		\label{FigSurgeryPrep}
	\end{center}
\end{figure}

\begin{Def}[adapted cut-off functions]
	\label{DefSurgeryCuttOffFunctions}
	In the situation of \Cref{ThmBaerDahlSurgeryPimped}, assume that for each $\varepsilon > 0$ (sufficiently small), we have a decomposition of $N$ into $N = U_{\varepsilon} \dot \cup A_{\varepsilon} \dot \cup V_{\varepsilon}$, where 	
	\begin{align}
		\label{EqDefUsrUsrrp}
		\begin{split}
			U_{\varepsilon} & :=  \{x \in N \mid \dist(x, S) < r_{_{\varepsilon}} \}, \\
			A_{\varepsilon} & := \{x \in N \mid r_{\varepsilon} \leq \dist(x,S) \leq r'_{\varepsilon} \},
		\end{split}
	\end{align}
	for some $r_{\varepsilon}, r'_{\varepsilon} > 0$. A family of cut-off functions $\chi^{\varepsilon} \in \mathcal{C}^\infty_c(N)$ is \emph{adapted} to these decompositions, if
	\begin{enumerate}
		\item 
			$0 \leq \chi^{\varepsilon} \leq 1$,
		\item
			$\chi^{\varepsilon} \equiv 0 \text{ on a neighborhood of } \bar U_{\varepsilon}$, 
		\item
			$\chi^{\varepsilon} \equiv 1 \text{ on } V_{\varepsilon}$, 
		\item
			$|\nabla \chi^{\varepsilon}| \leq \tfrac{c}{r_{\varepsilon}} \text{ on }N$ for some constant $c>0$.
	\end{enumerate}
	In case $\tilde N^{\varepsilon} = (N \setminus U_{\varepsilon}) \dot \cup \tilde U_{\varepsilon}$ is obtained from $N$ by surgery, the restriction $\chi^{\varepsilon}|_{N \setminus U_{\varepsilon}}$ can be extended smoothly by zero to a function $\chi^{\varepsilon} \in \mathcal{C}_c^\infty(\tilde N^{\varepsilon})$. The situation is depicted in \Cref{FigSurgeryPrep}.
\end{Def}

\begin{Rem}[cutting off spinor fields]
	\label{RemCuttingOffSpinorFields}
	We can use the cut-off functions from \Cref{DefSurgeryCuttOffFunctions} to transport spinor fields from $(N,g)$ to $(\tilde N^{\varepsilon}, \tilde g^{\varepsilon})$ and vice versa: For any $\psi \in L^2(\Sigma^g_{\K} N)$, we can think of $\chi^{\varepsilon} \psi$ as an element in $L^2(\Sigma^{\tilde g}_{\K} \tilde N^{\varepsilon})$ by extending $\chi^{\varepsilon} \psi$ to all of $\tilde N^{\varepsilon}$ by zero. Analogously, for any $\tilde \psi \in L^2(\Sigma^{\tilde g}_{\K} \tilde N^{\varepsilon})$, we can think of $\chi^{\varepsilon} \tilde \psi$ as an element in $L^2(\Sigma^g_{\K} N)$ by extending $\chi^{\varepsilon} \tilde \psi|_{N \setminus \tilde U_{\varepsilon}}$ by zero to $N$. This correspondence is not an isomorphism, but one does not loose ``too much'': It preserves smoothness and it is shown in \cite{BaerDahlSurgery} that under the assumptions of \Cref{ThmBaerDahlSurgeryPimped}, for each eigenspinor $\tilde \psi^{\varepsilon} \in L^2(\Sigma^{\tilde g^{\varepsilon}}_{\K} \tilde N^{\varepsilon})$ to an eigenvalue $\tilde \lambda^{\varepsilon} \in ]-\Lambda, \Lambda[$, the spinor field $\psi^{\varepsilon} := \chi^{\varepsilon} \tilde \psi^{\varepsilon} \in L^2(\Sigma^g_{\K} N)$ satisfies
	\begin{align}
		\label{EqBaerDahlEstimate}
		\| \Dirac^g_{\K} \psi^{\varepsilon}\|_{L^2(\Sigma^g_{\K} N)} &< (\Lambda + \tfrac{\varepsilon}{2}) \|\tilde \psi^{\varepsilon}\|_{L^2(\Sigma^{\tilde g}_{\K} \tilde N^{\varepsilon})}, \\
		\label{EqBaerDahlBelow}
		\|\psi^{\varepsilon}\|_{L^2(\Sigma^g_{\K} N)} & \geq \frac{\Lambda + \tfrac{\varepsilon}{2}}{\Lambda + \varepsilon}  \|\tilde \psi^{\varepsilon}\|_{L^2(\Sigma^{\tilde g}_{\K} \tilde N^{\varepsilon})}.
	\end{align}
	The proof of \cref{EqBaerDahlEstimate,EqBaerDahlBelow} is an integral part of the proof of \Cref{ThmBaerDahlSurgeryPimped}, see \cite[p. 69]{BaerDahlSurgery}. In combination, they imply the following crucial estimate for the Rayleigh quotient
	\begin{align}
		\label{EqDahlSurgeryRRBound}
		\frac{\| \Dirac^g_{\K} \psi^{\varepsilon}\|_{L^2(\Sigma^g_{\K} N)}^2}{\|\psi^{\varepsilon}\|_{L^2(\Sigma^g_{\K} N)}^2} < (\Lambda + \varepsilon)^2.
	\end{align}
	This estimate is then used to apply the \emph{min-max principle}, see \Cref{ThmMinMaxPrinciple}, which gives the conclusion of \Cref{ThmBaerDahlSurgeryPimped}. 
\end{Rem}

We will need the following version of \Cref{ThmBaerDahlSurgeryPimped} that is slightly more general. 

\begin{Thm}
	\label{ThmSurgeryFamilyPimped}
	Let $(N,\Theta)$ be a closed spin manifold of dimension $n \geq 3$, let $0 \leq l \leq n-3$ and $f$ be a surgery map with surgery sphere $S \subset N$ of dimension $l$ as in \Cref{DefNormalSurgery}. Let $(Z,\tau_Z)$ be a compact topological space, 
	\begin{align*}
		\mathbf{g}:(Z,\tau_Z) \to (\Rm(N), \mathcal{C}^2)
	\end{align*}
	be a continuous family of Riemannian metrics and let $\Lambda_1, \Lambda_2 \in \mathcal{C}^0(Z,\R)$, $\Lambda_1 < \Lambda_2$, such that
	\begin{align*}
		\forall z \in Z: \Lambda_1(z), \Lambda_2(z) \notin \spec \Dirac^{g_{z}}_{\K}.
	\end{align*}
	\begin{enumerate}
		\item
			For any $\varepsilon > 0$ (sufficiently small), there exists a spin manifold $(\tilde N^{\varepsilon}, \tilde \Theta^{\varepsilon})$, and a continuous family of Riemannian metrics
			\begin{align*}
				\mathbf{\tilde g^{\varepsilon}}:(Z, \tau_Z) \to (\Rm(\tilde N^{\varepsilon}), \mathcal{C}^2)
			\end{align*}
			such that for each $z \in Z$, the manifold $(\tilde N^{\varepsilon}, \tilde g^{\varepsilon}_{z})$ is obtained from $(N, g_{z})$ by surgery along $S$ and such that for all $z \in Z$, the operators $\Dirac^{g_{z}}_{\K}$ and $\Dirac^{\tilde g^{\varepsilon}_{z}}_{\K}$ are $(\Lambda_1(z), \Lambda_2(z), \varepsilon)$-spectral close. 
		\item
			For any open neighborhood $U \subset N$ of the surgery sphere, one can choose $\mathbf{\tilde g^{\varepsilon}}$ such that
			\begin{align}
				\label{EqMetricUnchangedSurgery}
				\forall z \in Z: \tilde g^{\varepsilon}_{z}|_{N \setminus U} = g_{z}|_{N \setminus U}.
			\end{align}
		\item
			If $\tilde \psi^{\varepsilon}_{z} \in L^2(\Sigma^{\tilde g^{\varepsilon}_{z}}_{\K} \tilde N^{\varepsilon})$, $z \in Z$, is in the span of the eigenspinors corresponding to eigenvalues in $[\Lambda_1(z), \Lambda_2(z)]$ and $\chi^{\varepsilon}_{z}$ is the cut-off function from \Cref{DefSurgeryCuttOffFunctions}, the spinor field $\psi^{\varepsilon}_{z} := \chi^{\varepsilon}_{z} \tilde \psi^{\varepsilon}_{z}$ satisfies
			\begin{align}
				\label{EqDahlSurgeryRRBoundPimped}
				\frac{\| (\Dirac^{g_{z}}_{\K}-c_z) \psi^{\varepsilon}_{z}\|_{L^2(\Sigma_{\K}^{g_{z}} N)}^2}{\|\psi^{\varepsilon}_{z}\|_{L^2(\Sigma^{g_{z}}_{\K} N)}^2} < (l_z + \varepsilon)^2,
			\end{align}		
			where $c_z:= \tfrac{1}{2}(\Lambda_1(z) + \Lambda_2(z))$, $l_z := \tfrac{1}{2}|\Lambda_2(z) - \Lambda_1(z)|$.		
	\end{enumerate}
\end{Thm}

\begin{figure}[t]
	\begin{center}
		\begin{tikzpicture}[scale=0.7]

\tikzset{
    partial ellipse/.style args={#1:#2:#3}{
        insert path={+ (#1:#3) arc (#1:#2:#3)}
    }
}

	\draw[thick] (6,0) circle (2);
	\draw[thick] (6,0) [partial ellipse=180:360:2 and 0.5];
	
	\draw[thick]
	(0,0)
	to [out=90, in=0] (-2.5,3)
	to [out=190, in=45] (-5,2.5)
	to [out=-135, in=90] (-6,0)
	to [out=-90, in=0] (-8.5,-1.5)
	to [out=180, in=135] (-8,-4.5)
	to [out=-45, in=-90] (-4.5,-2.5)
	to [out=90, in=180] (-3.5,-1.5)
	to [out=0, in=-90] (0,0);
	
	\draw[thick,rotate=-5] (-4.3,1) [partial ellipse=50:-68:0.5 and 1];
	\draw[thick,rotate=-6] (-3.8,0.8) [partial ellipse=115:232:0.5 and 1];
	
	\draw[thick,rotate=-30] (-4.3,-6) [partial ellipse=49:-58:0.5 and 1];
	\draw[thick,rotate=-30] (-3.7,-6.1) [partial ellipse=122:228:0.5 and 1];
	
	\fill[white] (-0.5,1) rectangle (0.5,0.5);	
	\fill[white] (4,1) rectangle (4.5,0.5);
		
	\draw[thick] (-1,1) to (5,1); 
	\draw[thick] (-1,0.5) to (5,0.5);	
	\draw[thick] [out=180,in=-90] (-1,1) to (-1.3,1.3); 
	\draw[thick] [out=180, in=90] (-1,0.5) to (-1.3,0.2); 
	\draw[thick] (5,1) [out=0, in=-90] to (5.3,1.3); 
	\draw[thick] (5,0.5) [out=0, in=90] to (5.3,0.2); 
	\draw[thick] (-1,0.75) ellipse (0.125 and 0.25); 
	\draw[thick] (5,0.75) ellipse (0.125 and 0.25);	
	
	\coordinate[label=above:$M \sharp S^m$] (A) at (2,-3);

\end{tikzpicture}
		\caption[Connected sum with a sphere]{Connected sum with a sphere.}
		\label{FigSphereSurgery}
	\end{center}
\end{figure}
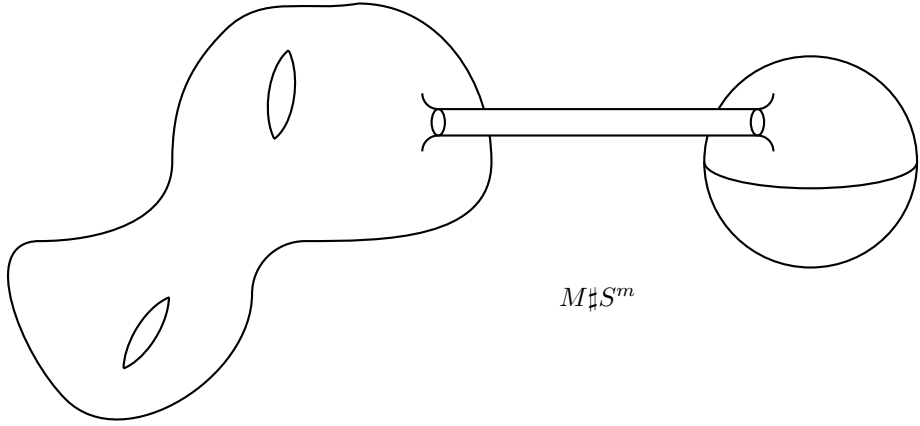

\begin{Rem}
	\Cref{ThmSurgeryFamilyPimped} generalizes \Cref{ThmBaerDahlSurgeryPimped} in the following ways.
	\begin{enumerate}
		\item 
			The metric $g$ is replaced by a compact $\mathcal{C}^2$-continuous family of metrics. It has already been observed by Dahl in a later paper, see \cite[Thm. 4]{DahlPresc}, that the proof of \Cref{ThmBaerDahlSurgeryPimped} goes through in this case.
		\item
			The interval $[-\Lambda, \Lambda]$ is replaced by the interval $[\Lambda_1, \Lambda_2]$, which might not be symmetric around zero. This is why one has to introduce $c$ and $l$ in \Cref{EqDahlSurgeryRRBoundPimped}.
		\item
			The field is $\K \in \{\R, \C\}$. This simply makes no difference in the proof.
		\item
			We replaced the constants $\Lambda_1$, $\Lambda_2$ by continuous functions on $Z$. This is possible, since their key function in the proof is to ensure that at any $z \in Z$ no eigenvalues enter or leave the spectral interval $[\Lambda_1(z), \Lambda_2(z)]$. Since they are continuous, they are also bounded, so uniform estimates are possible. (This generalization is not really needed in our proof of \Cref{MainThmHigher}.)
	\end{enumerate}
	These generalizations are all straightforward, but some more arguments can be found in \cite[A.8]{NiknoDiss}.
\end{Rem}

We are now able to prove the main result of this section.

\begin{Thm}[surgery stability]
	\label{ThmSurgeryStabilityOdd}
	Assume the hypothesis of \Cref{ThmSurgeryFamilyPimped} and in addition that $\mathbf{g}:Z \to Y(N)$, where $Y(N)$ is as in \Cref{DefPartEBundleGlobal} and let $E(N) \to Y(N)$ be the corresponding eigenbundle as in \Cref{ThmContEigenbundle}. Let $\mathbf{\tilde g}^{\varepsilon}:Z \to \Rm(\tilde N^{\varepsilon})$ be the loop from the conclusion of \Cref{ThmSurgeryFamilyPimped}. Then for $\varepsilon$ small enough, $\mathbf{\tilde g}^{\varepsilon}:Z \to Y(\tilde N^{\varepsilon})$ and the corresponding eigenbundle $E(\tilde N^{\varepsilon}) \to Y(\tilde N^{\varepsilon})$ satisfies
	\begin{align*}
		\mathbf{g}^*E \cong (\mathbf{\tilde g}^{\varepsilon})^* \tilde E^{\varepsilon}
	\end{align*}
	as vector bundles over $Z$. 
\end{Thm}

\begin{figure}[t] 
	\begin{center}
		\begin{tikzpicture}
		\tikzset{
			decoration={
				markings,
				mark=between positions 0.125 and 0.90 step 0.0625 with {
				\xdef\maxseq{\pgfkeysvalueof{/pgf/decoration/mark info/sequence number}}
				\coordinate (pt-\maxseq)
				at (0,0);},
			}
		}

		\draw[black,very thick] (0,0) ellipse (4 and 3);
	
		\draw[gray,very thick] (1,0.3) ellipse (4 and 3);
		
		\draw[gray] (0.8,0.1) ellipse (1.5 and 1);		

	\draw[black] (0.5,0) ellipse (1.5 and 1);
		
		\coordinate[label=above:{$Y$}] (A) at (-2.5,-2);
		\coordinate[label=above:{\color{gray}$\tilde Y$}] (A) at (4,1.5);
		\coordinate[label=above:{$X \setminus Y$}] (A) at (0.2,-0.7);
		\coordinate[label=above:{$\color{gray}\tilde X \setminus \tilde Y$}] (A) at (1,0);		
		\coordinate[label=above:{$\mathbf{g}$}] (A) at (-2.4,-0.4);
		\coordinate[label=above:{$\color{gray}\mathbf{\tilde g}$}] (A) at (-1.7,-1.4);				
		\coordinate[label=above:{$X=\Rm(N)$}] (A) at (-2,-3.6);
		\coordinate[label=above:{$\color{gray}\tilde X = \Rm(\tilde N)$}] (A) at (2,-3.6);
	
		\draw[thick,black,postaction={decorate}] (-2.4,-0.9) to (-2,-0.5)
			to [out=45,in=180] (0,-1.5)
			to [out=0,in=180] (2,-1)
			to [out=0,in=-90] (2.5,0)
			to [out=90,in=0] (1.5,1.6)
			to [out=170,in=0] (0.2,1.7)
			to [out=180,in=45] (-2,-0.2)
			to [out=220,in=45] (-2.4,-0.9);
		
		\foreach \c in {1,...,\maxseq}{
			\draw[black,rotate=1.35*360*(\c-1) / \maxseq] (pt-\c) +(0,-0.5) -- +(0,0.5);} 
				
		\fill[black] (-2.4,-0.9) circle (0.07);

		\draw[thick,gray,postaction={decorate}] (-2.1,-1) to (-1.9,-0.4)
			to [out=45,in=180] (0.1,-1.5)
			to [out=0,in=180] (2.1,-1.1)
			to [out=0,in=-90] (2.6,0)
			to [out=90,in=0] (1.4,1.5)
			to [out=170,in=0] (0,1.6)
			to [out=180,in=45] (-2,-0.1)
			to [out=220,in=45] (-2.1,-1);
		
	\foreach \c in {1,...,\maxseq}{
	\draw[gray,rotate=1.36*360*(\c-1) / \maxseq] (pt-\c) +(0,-0.5) -- +(0,0.5);
	}
	
	\fill[gray] (-2.1,-1) circle (0.07);
			
\end{tikzpicture}
		\caption[Surgery Stability]{Lasso before surgery and {\color{gray}after surgery}}
		\label{FigSurgStable}
	\end{center}
\end{figure}
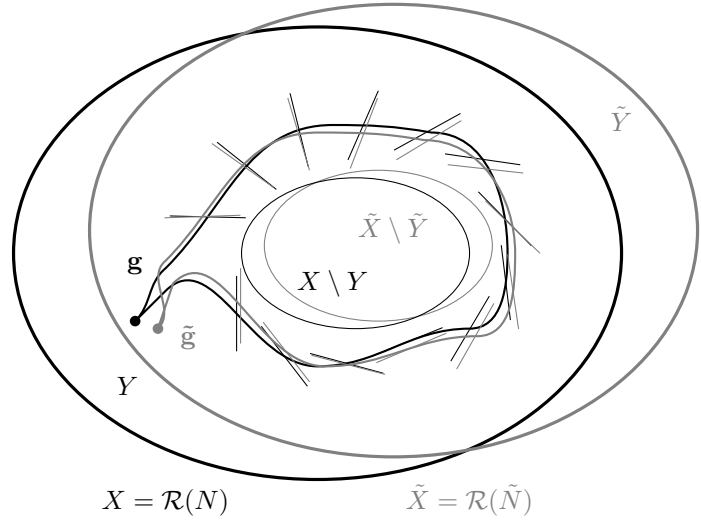

\begin{Prf}
	The idea of the proof is to show that the situation of the Lasso Lemma (\cref{LemLasso}) still holds, if one perturbs everything a bit, see \Cref{FigSurgStable}.	
	\begin{steplist}
		\step[$\mathbf{\tilde g}^{\varepsilon}:Z \to Y(\tilde N^{\varepsilon})$]
			By the conclusion of \Cref{ThmSurgeryFamilyPimped}, we obtain that for each $z \in Z$, the operators $\Dirac^{g_{z}}_{\K}$ and $\Dirac^{\tilde g_{z}^{\varepsilon}}_{\K}$ are $(\Lambda_1(z), \Lambda_2(z), \varepsilon)$-spectral close. We fix some $z_0 \in Z$ and take an enumeration of $\Dirac^{\tilde g^{\varepsilon}_{z_0}}_{\K}$ such that $\tilde \lambda_1^{\varepsilon}(\tilde g^{\varepsilon}_{z_0})$ is the smallest eigenvalue $ > \Lambda_1(z_0)$. Let $\{\tilde \lambda_j^{\varepsilon}\}_{j \in \Z}$ be the corresponding enumeration of eigenvalues on $\tilde N^{\varepsilon}$ as in \Cref{ThmEvpaper}. We obtain, that for $\varepsilon$ small enough
			\begin{align*}
				\forall z \in Z:
				\tilde{\lambda}^{\varepsilon}_{0}(\tilde g^{\varepsilon}_{z})
				\leq \Lambda_1(z) 
				< \tilde \lambda^{\varepsilon}_1(\tilde g^{\varepsilon}_{z}) 
				\leq \ldots 
				\leq \tilde \lambda^{\varepsilon}_k(\tilde g^{\varepsilon}_{z}) 
				< \Lambda_2(z)
				\leq \tilde \lambda^{\varepsilon}_{k+1}(\tilde g^{\varepsilon}_{z}).
			\end{align*}
			Therefore, for each $z \in Z$, there are $k$ eigenvalues between $\Lambda_1(z)$ and $\Lambda_2(z)$. Since $k$ is odd, at least one eigenvalue must have odd multiplicity. Hence $\mathbf{\tilde g}^{\varepsilon}:Z \to \tilde Y^{\varepsilon}(N)$.
						
		\step[passing from $\tilde N^{\varepsilon}$ to $N$] 
			Let $L^2(\Sigma_{\K} N) \to \Rm(N)$ be the universal spinor field bundle, see \Cref{ThmUniSpinorFieldBdle}, $E(N)$ be the bundle from \Cref{DefPartEBundleGlobal} and define
			\begin{align*}
				\mathcal{H} &:= \mathbf{g}^*(L^2(\Sigma_{\K} N)) \to Z, \\
				E &:=  \mathbf{g}^*(E(N)) \to Z, \\
				\tilde E^{\varepsilon} &:= (\mathbf{\tilde g^{\varepsilon}})^*(E(\tilde N^{\varepsilon})) \to Z.
			\end{align*}
			For any $z \in Z$, let $\chi^{\varepsilon}_z$, be the canonical cut-off functions from \Cref{DefSurgeryCuttOffFunctions}. These functions can be chosen such that $\chi^{\varepsilon}_z$ depends continuously on $z$. We consider the map
			\begin{align*}
				P:\tilde E^{\varepsilon} \to \mathcal{H}, && 
				\tilde \psi_z^{\varepsilon} \in \Gamma(\Sigma^{\tilde g_z^{\varepsilon}}_{\K} \tilde N^{\varepsilon}) \mapsto \psi_z^{\varepsilon} := \chi^{\varepsilon}_z \tilde \psi_z^{\varepsilon} \in \Gamma(\Sigma^{g_z}_{\K} N),				
			\end{align*}
			see also \Cref{RemCuttingOffSpinorFields}. In case $P(\tilde \psi^{\varepsilon}_z) = 0$, we obtain $\tilde \psi^{\varepsilon}_z = 0$ by the weak unique continuation property for Dirac type operators, see \cite[Rem. 2.3c)]{BoossWeakUCP}. Therefore, this is a continuous morphism of vector bundles that is fibrewise injective. Hence, $P$ has constant rank and its image $E^{\varepsilon} := P(\tilde E^{\varepsilon}) \subset \mathcal{H} \to Z$ is a continuous vector bundle isomorphic to $\tilde E^{\varepsilon}$. 
			
		\step[analyze Rayleigh quotients]
			We consider the $\lambda_j$'s as functions on $Z$ by pulling them back via $\mathbf{g}$ (and analogously for $\tilde \lambda_j^{\varepsilon}$). For any $z \in Z$, if $\lambda_1(z), \ldots, \lambda_k(z)$ are the eigenvalues of $\Dirac^{g_z}_{\K}$ in $[\Lambda_1(z), \Lambda_2(z)]$, then $(\lambda_1(z)-c_z)^2, \ldots, (\lambda_k(z) - c_z)^2)$ are the eigenvalues of $(\Dirac^{g_z}_{\K}-c_z)^2$ in $[0,l(z)^2]$, where $c_z := \tfrac{1}{2} \left( \Lambda_1(z) + \Lambda_2(z) \right)$ and $l_z := \tfrac{1}{2} |\Lambda_2(z) - \Lambda_1(z)|$. The span of their collective eigenspinors is the same space $E_z$. It follows from \Cref{EqDahlSurgeryRRBoundPimped} that the Rayleigh quotients satisfy
			\begin{align*}
				\forall z \in Z: \forall \psi^{\varepsilon}_z \in \Gamma (\Sigma^{g_z}_{\K} N): 
				\frac{\| (\Dirac^{g_z}_{\K} -c_z) \psi^{\varepsilon}_z\|_{L^2(\Sigma^{g_z}_{\K} N)}^2 }{\|\psi^{\varepsilon}_z\|_{L^2(\Sigma^{g_z}_{\K} N)}^2} < (l_z + \varepsilon)^2.
			\end{align*}
			Now, choose $\varepsilon$ small enough such that for all $z \in Z$, we have $l_z + \varepsilon < \rho_{k+1}(z)$, where $\rho_{k+1}(z)$ is the $(k+1)$-th eigenvalue of $(\Dirac^{g_z}_{\K} - c_z)^2$. This is possible due to the continuity of $\mathbf{g}$ and since
			\begin{align*}
				\forall z \in Z: \lambda_0(z) \leq \Lambda_1(z) < \lambda_1(z) \leq \ldots \leq \lambda_k(z) < \Lambda_2(z) \leq \lambda_{k+1}(z)
			\end{align*}
			by hypothesis. By \Cref{ThmRayleighRitzDistance}, we obtain
			\begin{align*}
				\forall z \in Z: \forall 1 \leq j \leq k: d(E_z,\psi^{\varepsilon}_z)^2 \leq \frac{l(z) + \varepsilon}{\rho_{k+1}(z)} < 1.
			\end{align*}
			By \Cref{ThmSignStability}, we obtain that $E \cong E^{\varepsilon} \cong \tilde E^{\varepsilon}$.
	\end{steplist}
	All in all, we have verified the hypothesis of \Cref{LemLasso} for $E(\tilde N^{\varepsilon}) \to Y(\tilde N^{\varepsilon})$ and $\mathbf{\tilde g}^{\varepsilon}$ , which proves the claim.
\end{Prf}

\begin{Rem}
	\Cref{ThmSurgeryStabilityOdd} holds also for other sets than $Y(N)$. For instance it holds for the larger set
	\begin{align*}
		\{g \in \Rm(N) \mid \lambda_{0}(g) < \lambda_1(g), \lambda_{k}(g) < \lambda_{k+1}(g) \},
	\end{align*}	
	with essentially the same proof. It does not hold for arbitrary subsets of $\Rm(N)$ though. For instance, if one replaces the condition $\exists 1 \leq j \leq k: \lambda_j \in 2 \N+1$ by $\exists 1 \leq j \leq k: \lambda_j \in 2 \N$, the first step of the proof no longer holds, because an eigenvalue of even multiplicity might split up into two eigenvalues of odd multiplicity. 
\end{Rem}

\section{Proof of \Cref{MainThmHigher}}

We are now in a position to put all the results together.

\begin{Prf}[of \Cref{MainThmHigher}]
	The idea is to apply the surgery stability theorem \Cref{ThmSurgeryStabilityOdd}, to the connected sum $M \sharp \S^m$, see \Cref{FigSphereSurgery}.
	\begin{steplist}
		\step[build lasso on sphere]
			By \Cref{ThmRotationsOddLoop}, there exists an odd loop of spin diffeomorphisms on $\S^m$. By \Cref{ThmNbhdRoundSphereOdd}, there exists a metric on $\S^m$, for which at least one eigenvalue $\lambda$ as an odd multiplicity $k$. The associated loop of metrics $\mathbf{g}_{\S^m}$ is a $\mathcal{C}^2$-continuous map $\mathbf{g}_{\S^m}:\S^1 \to Y(\S^m)$, where $Y(\S^m)$ is as in \Cref{DefYKM}. We obtain the associated real eigenbundle $E(\S^m)$ as in \Cref{DefPartEBundleGlobal}. By \Cref{ThmContEigenbundle}, $E(\S^m) \to Y(\S^m)$ is a continuous vector bundle of rank $k$. By \Cref{ThmS1ActionBundleTwisted}, this bundle it not trivial. This verifies the hypothesis of \Cref{LemLasso} on the sphere.
		\step[prepare $M$ for surgery]
			Define the manifold $N := M \amalg \S^m$. Choose any metric $\tilde g$ on $M$. We obtain the loop $\mathbf{g} := \tilde g \amalg \mathbf{g}_{\S^m}:\S^1 \to \Rm(N)$. In case $\lambda \in \spec \tilde g$, we first scale $\mathbf{g}_{\S^m}$ a little such that $\lambda \notin \spec \tilde g$. We obtain that $\spec g_z$ is constant with respect to $z \in \S^1$. Since the spectrum is also discrete, we can certainly find continuous (even constant) functions $\Lambda_1, \Lambda_2: \S^1 \to \R$ such that
			\begin{align*}
				\forall z \in \S^1: \lambda_0(g_z) \leq \Lambda_1 < \underbrace{\lambda_1(g_z) = \ldots = \lambda_k(g_z)}_{=\lambda} < \Lambda_2 \leq \lambda_{k+1}(g_z).
			\end{align*}
			It follows that $\mathbf{g}:\S^1 \to Y(N)$. We get that $E(N) \to Y(N)$ is a non-trivial vector bundle of real rank $k$, since it is isomorphic to the pullback of $E(M) \to Y(M)$ along the map $\Rm(N) \to \Rm(M)$, $g \amalg h \mapsto g$. 
			
		\step[perform surgery]
			We extend the loop $\mathbf{g}:\S^1 \to Y(N)$ to a disc $\D^2 \to \Rm(N)$, which is still denoted by $\mathbf{g}$, and set $Z := \D^2$. We apply \Cref{ThmSurgeryFamilyPimped} to $N$ in dimension $l=0$ for $\K = \R$, i.e. we obtain a connected sum $\tilde N^{\varepsilon} = M \sharp \S^m$ together with a resulting family $\mathbf{\tilde g}^{\varepsilon}:\D^2 \to \Rm(\tilde N^{\varepsilon})$ of metrics. By \Cref{ThmSurgeryStabilityOdd}, this gives a loop $\mathbf{\tilde g}^{\varepsilon}|_{\S^1}:\S^1 \to Y(\tilde N^{\varepsilon})$ and the corresponding bundle $E(\tilde N^{\varepsilon}) \to Y(\tilde N^{\varepsilon})$ satisfies $(\mathbf{\tilde g}^{\varepsilon}|_{\S^1})^* E(\tilde N^{\varepsilon}) \cong (\mathbf{g}|_{\S^1})^*(E(N))$, thus it is also not trivial.
	
	\end{steplist}
	All in all, we have verified the hypothesis of \Cref{LemLasso} on $\tilde N^{\varepsilon} = M \sharp \S^m$, which is spin diffeomorphic to $M$. This proves the first part of \Cref{MainThmHigher}. Recall from \Cref{RemDiscExtension} that our metric lies somewhere on the disc $\mathbf{\tilde g}^{\varepsilon}:\D^2 \to \Rm(M)$. Therefore, the second claim follows from \Cref{EqMetricUnchangedSurgery}.
\end{Prf}

\begin{appendix}

\section{Rayleigh Quotients and the min-max principle}
\label{SectRayleighRitz}

In this section, we provide a version of the min-max principle suitable for our needs, see also \cite[XIII.1]{Reed4}.

\begin{Thm}[min-max principle]
	\index{min-max principle}
	\label{ThmMinMaxPrinciple}
	Let $H$ be a Hilbert space and $T: \dom(T) \subset H \to H$ be a densely defined self-adjoint operator with compact resolvent. Assume that the spectrum of $T$ satisfies $b \leq \lambda_1 \leq \lambda_2 \leq \ldots$ for some lower bound $b \in \R$. Then for any $k \in \N$
	\begin{align}
		\label{EqMinMaxEV}
		\lambda_k = \min_{\substack{U \subset \dom(T), \\ \dim(U) = k}}{\max_{\substack{x \in U, \\ \|x\|=1}}{\langle Tx, x \rangle }},
	\end{align}
	where the $\min$ is taken over all linear subspaces $U \subset \dom(T)$ of dimension $k$.
\end{Thm}

\begin{Thm}
	\label{ThmRayleighRitzDistance}
	Let $H$ be a Hilbert space, $T:H \to H$ be a densely defined operator, self-adjoint with compact resolvent. Let $\Lambda> 0$, $k \in \N$ and assume that the first $k+1$ distinct eigenvalues of $T$ satisfy
	\begin{align*}
		0 \leq \lambda_1 < \ldots < \lambda_k < \Lambda < \lambda_{k+1}.
	\end{align*}
	Define $E^{(\nu)} := \ker(T-\lambda_{\nu})$, $V := \bigoplus_{\nu=1}^{k}{E^{(\nu)}}$, and let $x \in \dom(T)$, $\|x\|=1$, such that $\langle Tx, x \rangle  \leq \Lambda + \varepsilon$. Then the distance between $x$ and $V$ satisfies $\dist(V,x)^2 \leq \frac{\Lambda + \varepsilon}{\lambda_{k+1}}$.
\end{Thm}

\begin{Prf}
	Consider the orthogonal decomposition
	\begin{align*}
		x = \sum_{\nu=1}^{\infty}{x^{(\nu)}}, && x^{(\nu)} \in E^{(\nu)}.
	\end{align*}
	By hypothesis
	\begin{align*}
		\varepsilon + \Lambda & \geq \langle Tx, x \rangle 
		=\langle \sum_{\nu=1}^{k}{Tx^{(\nu)}} + \sum_{\nu=k+1}^{\infty}{Tx^{(\nu)}}, \sum_{\nu=1}^{k}{x^{(\nu)}} + \sum_{\nu=k+1}^{\infty}{x^{(\nu)}} \rangle \\
		&=\langle \sum_{\nu=1}^{k}{\lambda_{\nu} x^{(\nu)}} + \sum_{\nu=k+1}^{\infty}{\lambda_{\nu} x^{(\nu)}}, \sum_{\nu=1}^{k}{x^{(\nu)}} + \sum_{\nu=k+1}^{\infty}{x^{(\nu)}} \rangle \\	
		&=\sum_{\nu=1}^{k}{\lambda_{\nu} \|x^{(\nu)}\|^2} + \sum_{\nu=k+1}^{\infty}{\lambda_{\nu} \|x^{(\nu)}\|^2} \\
		& \geq \lambda_{k+1} \sum_{\nu=k+1}^{\infty}{\|x^{(\nu)}\|^2}.
	\end{align*}
	Let $P_V:H \to H$ be the orthogonal projection onto $V$. We obtain
	\begin{align*}
		\dist(V, x)^2
		= \|P_V(x) - x\|^2 
		=\sum_{\nu=k+1}^{\infty}{\|x^{(\nu)}\|^2}
		\leq \frac{\Lambda + \varepsilon}{\lambda_{k+1}}.
	\end{align*}
\end{Prf}

\end{appendix}	
		
		\setlength{\nomitemsep}{-0.8\parsep}   
		\printnomenclature[8.5em]
			
		\printbibliography[title=References]
		\addcontentsline{toc}{section}{References}
		
\end{document}